\input ./style/arxiv-general.cfg
\documentclass[aos,MSNbibl,seceqn,nameyear,dvips]{arximspdf}
\makeatletter
   \@ifpackageloaded{graphicx}{}{\usepackage{graphicx}}
\makeatother
\usepackage{dcolumn}

\doi{10.1214/15-AOS1378}% Updated by VTEXPTS2LaTeX.exe, 21.10.2015
\volume{44}
\issue{2}
\pubyear{2016}
\firstpage{598}
\lastpage{628}
\docsubty{FLA}

\makeatletter
\newcolumntype{d}[1]{D{.}{.}{#1}}
\def\implies{\to}

\newproclaim{rem}{Remark}

\newproclaim{remark}{Remark}[section]
\newtheorem{theorem}{Theorem}[section]

\newtheorem{lemma}{Lemma}[section]

\newtheorem{cor}{Corollary}[section]
\newproclaim{definition}{Definition}[section]
\newproclaim{Example}{Example}[section]

\newproclaim{example}{Example}

\newproclaim{assumption}{Assumption}[section]

\newcommand{\Tr}{\operatorname{Tr}}

\newproclaim{defn}{Definition}
\newproclaim{Remark}{Remark}
\makeatother

\begin{document}
\begin{frontmatter}

\title{Large sample behaviour of high dimensional autocovariance matrices}
\runtitle{High dimensional autocovariance matrices}

\begin{aug}
% Corresponding author: Arup Bose - bosearu@gmail.com% Updated by
%VTEXPTS2LaTeX.exe, 21.10.2015 08:17
\author[A]{\fnms{Monika} \snm{Bhattacharjee}\ead[label=e1]{monaiidexp.gamma@gmail.com}}
\and
\author[A]{\fnms{Arup} \snm{Bose}\corref{}\thanksref{T2}\ead[label=e2]{bosearu@gmail.com}}
\runauthor{M. Bhattacharjee and A. Bose}
\thankstext{T2}{Supported by J.~C. Bose National Fellowship,
Department of Science and Technology, Government of India.}
\affiliation{Indian Statistical Institute}
%\dedicated{}
\address[A]{Statistics and Mathematics Unit\\
Indian Statistical Institute\\
203 B.T. Road, Kolkata 700108\\
India\\
\printead{e1}\\
\phantom{E-mail:\ }\printead*{e2}}
\end{aug}

% HISTORY:
%
\received{\smonth{6} \syear{2015}}% Updated by VTEXPTS2LaTeX.exe,
%21.10.2015 08:17
%
\revised{\smonth{8} \syear{2015}}% Updated by VTEXPTS2LaTeX.exe,
%21.10.2015 08:17

% ABSTRACT
%
\begin{abstract}
The existence of limiting spectral distribution (LSD) of $\hat{\Gamma
}_{u} + \hat{\Gamma}_{u}^{*}$, the symmetric sum of the sample
autocovariance matrix $\hat{\Gamma}_{u}$ of order $u$, %for every $i$
is known
when the observations are from an infinite dimensional vector linear
process with appropriate (strong) assumptions on the coefficient matrices.
Under significantly weaker conditions,
we prove, in a unified way, that %the expected average trace of any
%polynomial in these matrices converge.
%In particular,
the LSD of any symmetric polynomial in these matrices such as
$\hat{\Gamma}_{u} + \hat{\Gamma}_{u}^{*}$, $\hat{\Gamma}_{u}\hat
{\Gamma
}_{u}^{*}$,
$\hat{\Gamma}_{u}\hat{\Gamma}_{u}^{*}+\hat{\Gamma}_{k}\hat
{\Gamma
}_{k}^{*}$ exist.
Our approach is through the more intuitive algebraic method of free
probability in conjunction with the method of moments. Thus, we are
able to provide a general description for the limits in terms of some
freely independent variables.
All the previous results follow as special cases.
We suggest statistical uses of these LSD and related results in order
determination and white noise testing.
\end{abstract}

% KEYWORDS
\end{frontmatter}

%s1 #&#
\section{Introduction} \label{sec: Intro}

Multivariate linear time series models such as the
Autoregressive Moving Average (ARMA) processes are fundamental in the
theory of econometrics and finance. Moreover, time series data where
the dimension grows along with the sample size are becoming
increasingly frequent.
%Spatial data over time, fMRI data, spectroscopic imaging are some
%examples of such high-dimensional time %series data. %where such
%models are used.
A key model in these situations is %the
%infinite dimensional vector autoregressive (IVAR) and the
%infinite dimensional vector autoregressive moving average (IVARMA)
%processes or more generally,
%. These are special cases
%(see \citep{BB2013a}) of the more general
the infinite dimensional
moving average process of order infinity,
$\operatorname{MA}(\infty)$, where the sample $\{ X_{t.p}^{(n)}: t=1,2,\ldots, n\}$
of size
$n$ satisfies
%
%e1.1 #&#
\begin{equation}
\label{eqn:model} X_{t.p}^{(n)} =\sum
_{j=0}^\infty\psi_{j.p}^{(n)}
\varepsilon_{t-j.p}\qquad \forall t, n \ge1 \mbox{ (almost surely)}.
\end{equation}

For all $t$, $X_{t.p}^{(n)}$ and $\varepsilon_{t.p}$ are
$p$-dimensional vectors
and $\psi_{j.p}^{(n)}$ are $p\times p$ \textit{coefficient matrices}
and $\psi_{0.p}^{(n)}=I_{p}$. Precise assumptions of independence,
finiteness of moments and conditions on the matrices are discussed later.

We work in the framework of the particular
high dimensional model where the dimension $p$ increases
proportionately with the sample size $n$, so that $p=p(n) \rightarrow
\infty$ and $\frac{p}{n} \rightarrow y \in(0,\infty)$. The infinite
sum in (\ref{eqn:model}) exists in the almost sure sense under suitable
decay conditions on $\{\psi_{j.p}^{(n)}\}$. If $\psi_{j.p}^{(n)}=0, \forall j>q$, then it will be called an
$\operatorname{MA}(q)$ process. For convenience, we will write $p$ for $p(n)$ and
$\psi
_j, \varepsilon_t$ and $X_t$, respectively, for $\psi_{j.p}^{(n)},
\varepsilon_{t.p}$ and $X_{t.p}^{(n)}$. %Initially we assume that %on
%which we impose the following
%assumption.
Many researchers have worked on this model recently. See, for example,
\citeauthor{FHLR2000} (\citeyear{FHLR2000,FHLR2004}), \citet{FL2001}, \citet{BB2013a},
\citet
{JWBN2014} and \citet{LAP2013}.
%\vskip5pt
%$\textbf{(A3)}$ For all $i,t \ge1$, all the moments of $
%\varepsilon_{t,i}$'s
%exist.
%\vskip5pt
%Later we will get rid of assumption (A3) for some special cases (see
%Remark~\ref{rem: truncation}). %and show that the existence of lower
%order moments are actually enough.
%Also, all the discussions below will not be changed too much if the
%dispersion matrix of $\varepsilon_{t,p}$ is any positive definite
%matrix $\Sigma_{p}$ instead of $I_p$. In this context, the only change
%is to replace the coefficient matrices $\psi_j$ by $\psi_j
%\Sigma^{1/2}_{p}$ for all $j=0,1,2,\ldots.$
%\textbf{Problem 1: Truncation may reduce this moment
%condition.}

One of the key quantities in time series analysis is the autocovariance
matrix. % and their estimates.
The \textit{population autocovariance matrices}
are defined as
\[
\Gamma_{u.p}: = E \bigl(X_{t.p} X^{*}_{(t+u).p}
\bigr) = \sum_{j=1}^\infty\psi_j
\psi_{j+u}^{*},\qquad u=0, 1, \ldots.
\]
%
%Here our goal is to find out the
%limiting spectral distribution of $\hat{\Gamma}_0$,
%$\hat{\Gamma_i} \hat{\Gamma_i}^\prime$ and $(\hat{\Gamma}_i +
%\hat{\Gamma}_i^\prime)$ for all $i=1,2,3 \ldots$

The moment estimator of $\Gamma_{u.p}$ is the \textit{sample
autocovariance matrix},
%
%e1.2 #&#
\begin{equation}
\label{eqn: sampleautocovdefn} \hat{\Gamma}_{u.p} = \frac{1}{n} \sum
_{t=1}^{n-u} X_{t.p} X_{(t+u).p}^{*},\qquad
0 \le u \le n-1.
\end{equation}
We often write
$\Gamma_u$ and $\hat{\Gamma}_{u}$, respectively, for $\Gamma_{u.p}$
and $\hat{\Gamma}_{u.p}$. Our goal is to study the large sample
behaviour of the random matrices $\hat{\Gamma}_{u}$ and use the
asymptotic results for statistical inference purposes such as order
determination of infinite dimensional moving average or autoregressive
processes. Since we are dealing with several matrices of
increasing dimension together, we need to give precise meaning to the
large sample behaviour in our context.

The most common way to capture the large sample behaviour of a sequence
of random matrices is through its spectral distribution.
The empirical spectral distribution (ESD) of an $n \times n$ (random)
matrix $R_n$ is the
(random) probability distribution with mass $1/n$ at each of its
eigenvalues. If it
converges weakly (almost surely) to a
(non-degenerate) probability distribution, then the latter is called
the \textit{limiting spectral distribution}
(LSD) of $R_n$.
%Suppose $R_n$ is an $n \times n$ real symmetric matrix. Let $\lambda_1,
%\lambda_2,\ldots, \lambda_n \in\R$ denote its eigenvalues. The
%\textit{spectral measure} $\mu_n$ of $R_n$ is the measure on $\R$
%given by
%\begin{equation} \label{eqn: ESD} \nonumber
%\mu_n = \frac{1}{n} \sum_{i=1}^n \delta_{\lambda_i}
%\end{equation}
%where $\delta_x$ is the Dirac delta measure at $x$. The
%probability distribution function on $\R$ corresponding to
%$\mu_n$, is known as the \textit{Empirical Spectral Distribution}
%(ESD) of
%$R_n$. If this ESD converges weakly (either almost surely or in
%probability) to some unique
%(non-degenerate) probability distribution, then the limit is called
%the \textit{limiting spectral distribution} (LSD)
%of $R_n$.
Incidentally, the study of the limit spectrum of non-Hermitian matrices
is extremely difficult and very few results are known for general
non-Hermitian sequences. Researchers have concentrated on the additive
symmetrized version $\hat{\Gamma_u} +\hat{\Gamma}_u^*$, %of the
%autocovariance matrices,
one at a time.
%Later we shall introduce methods from free probability to consider, in
%some sense, the joint behaviour of these matrices. %We do not deal
%with the spectrum of the non-symmetric autocovariance matrices in this
%article but with symmetric matrix polynomials in them.

%Though this needs further investigation to be formalized,
%nevertheless, it
%This provides another instance of the importance of the LSD of $\{\hat{
%\Gamma}_{k}\}$ in inferential problems.

%Thus a detailed study of the LSD of $\{\hat{\Gamma}_{k}\}$ has
%significant statistical importance in high-dimensional time series.
%Our goal is to provide the large sample behaviour of the
%autocovariance matrices.

One widely used approach to establish the LSD
%of a real symmetric random matrix is %the method of moments, as
%outlined in Lemma~\ref{lem:LSDmoment}.
%For some applications of this method, see
%one may refer to %\citet{AGZ2009},
% \citet{B1999} and \citet{BBS2013}.
%%\citet{BS2009}.
%Limiting Spectral Distribution of any matrix is the weak limit of the
%Empirical Spectral Distribution of the same matrix defined in (
%\ref{eqn: ESD}).
%There is a huge random matrix literature, where the LSD of different
%type of random matrices are discussed.
%There are two main techniques that have been widely used in the random
%matrix theory for establishing LSD of random matrices --- the method
%of moments and the method of Stieltjes transforms. A brief description
%of these two methods are given below.
%\vskip5pt
%\textbf{Moment method:}
%\textbf{Method based on Stieltjes transformation:}
%Another widely used method is
is that of Stieltjes transformation, which for any finite measure $\mu$
on the real line equals
%
%e1.3 #&#
\begin{equation}
\label{eqn: stieltjes} m_{\mu}(z) = \int\frac{1}{x-z} \mu(dx),\qquad  z \in\mathbb
{C}^{+}:= \{ x+iy: x \in\mathbb{R}, y>0\}.
\end{equation}
%
%where $\mathbb{C}^{+} $. %This transformation has one-to-one
%correspondence with the measures on real line.
Pointwise convergence of Stieltjes transforms to a Stieltjes transform
implies the convergence of the corresponding distributions. In random
matrix theory, this convergence is proved by linking the Stieltjes
transform of the ESD to the resolvent and showing convergence by
martingale convergence methods.
See, for example, \citet{BS1995}, \citet{S1995}, \citet{BZ2008} and
Bai and Silverstein (\citeyear{BS2009}). % All the existing works on high-dimensional time
%series were based on the LSD of sample covariance matrices.
All existing works regarding LSD of autocovariance matrices are based
on this method.

Let us first discuss briefly these results.
Consider the simplest case of (\ref{eqn:model}) where $X_t =
\varepsilon
_t = (\varepsilon_{t,1},\varepsilon_{t,2},\ldots,\varepsilon
_{t,p})^{T}$ and $\{\varepsilon_{t,j}\}$ are i.i.d. with mean zero and
variance one and with enough high moments.
%Then under appropriate assumptions on moments of $\{\varepsilon_t\}$,
%following two results are known in the literature. \vskip5pt
Then $\hat{\Gamma}_{0}$ is nothing but the unadjusted sample
variance--covariance matrix, and it is well known that its LSD is the
Mar\v{c}enko--Pastur law. See, for example, \citet{MP1967} and \citet{BS2009}.
For the same model, \citet{JWBN2014} showed that the LSD of $\hat
{\Gamma
}_{u} + \hat{\Gamma}_{u}^{*}$ exist and are equal for every $u \geq1$
and derived its Stieltjes transform. % when $X_t = \varepsilon_t\
%\forall t$ satisfies appropriate conditions. More details are
%discussed later in ????

\citet{PS2011} and \citet{Y2012} derived the limiting Stieltjes
transformation of $\{\hat{\Gamma}_{u} + \hat{\Gamma}_{u}^{*}\}$, $u
\geq0$, when the components of $X_{t}$ are independent samples from an
identical univariate $\operatorname{MA}(q)$ process. \citet{LAP2013} appears to be the
only work in model (\ref{eqn:model}) for arbitrary $q$. For each $u$,
they established the existence of the LSD of $\{\hat{\Gamma}_{u} +
\hat
{\Gamma}_{u}^{*}\}$ and derived its Stieltjes transform as a solution
of a pair of functional equations. To derive this result, they assumed
that $\{\varepsilon_{t,j}\}$ are i.i.d. with finite $4$th moment.

Their assumptions on $\{\psi_j\}$ are, however, quite restrictive. They
assumed that
$\{\psi_j\}$ are Hermitian and simultaneously diagonalizable (the
latter assumption can be replaced by the assumption that $\{\psi_j\}$
are Toeplitz matrices with suitable decay conditions on their entries).
Even so, this excludes many interesting linear processes (such as model
$4$ in Section~\ref{subsec: example}). %As a simple example, consider
%the $\operatorname{MA}(2)$ process
%\begin{equation} \label{eqn: ma1example}
%X_{t,p} = \varepsilon_{t,p} + A_{p}\varepsilon_{t-1,p}+B_{p}
%\varepsilon_{t-2,p},
%\end{equation}
% where $A_{p} = ((a_{i,j}))_{p \times p}$ with
%$a_{i,n+1-i}=1$ and $0$ otherwise and $B_{p} = ((b_{i,j}))_{p \times
%p}$ is diagonal with $b_{i,i}=I(1\leq i \leq[\frac{n}{2}])-I([
%\frac{n}{2}]< i \leq n)$.
%This model is excluded
%. It is easy to see that $\psi_{0,p} = I_p$ and $A_{p}$ are not
%simultaneously diagonalizable and also spatial dependence is present
%in the model.
%under the set up of \citet{LAP2013} since $A_p$ and $B_p$ are neither
%Toeplitz nor simultaneoulsy diagonalizable. %Hence their result cannot
%be used to obtain the LSD of $\{\hat{\Gamma}_{k}+\hat{\Gamma}_{k}^{*}
%\}$ in model (\ref{eqn: ma1example}).
%Our result on the LSD (described later) is able to include this model.
To indicate another limitation of this assumption,
%of \citet{LAP2013},
suppose further that $\varepsilon_{t} \sim\mathcal{N}(0,I)$.
%and $\{\psi_j\}$ are hermitian and simultaneously diagonalizable.
Let $U$ be a unitary matrix such that $U\psi_j U^{*} =: \Lambda_j$
(say) are diagonal matrices.
%Let $U$ be a unitary matrix such that $\{U\psi_j U^{*} =: \Lambda_j\}$
%(say) are diagonal matrices.
Since $U\varepsilon_{t}$ and $\varepsilon_{t}$ are identically
distributed and $UU^{*}=I$, as far as the LSD of $\hat{\Gamma
}_{u}+\hat
{\Gamma}_{u}^{*}$ is concerned, (\ref{eqn:model}) is equivalent to
the model
%\begin{eqnarray} \label{eqn: diagassump}
%&& n (\hat{\Gamma}_{k}+\hat{\Gamma}_{k}^{*}) = \sum X_t X_{t-k}^{*} +
%\sum X_{t-k} X_{t}^{*} \nonumber\\
%&=& \sum_{t,j_1,j_2} \psi_{j_1} \varepsilon_{t-j_1}
%\varepsilon_{t-j_2-k}^{*}\psi_{j_2} + \sum_{t,j_1,j_2} \psi_{j_2}
%\varepsilon_{t-j_2-k} \varepsilon_{t-j_1}^{*}\psi_{j_1}\nonumber\\
%& =& U^{*} \bigg( \sum_{t,j_1,j_2} (U\psi_{j_1}U^{*}) U
%\varepsilon_{t-j_1} \varepsilon_{t-j_2-k}^{*}U^{*}(U\psi_{j_2}U^{*}) +
%\nonumber\\
%&& \sum_{t,j_1,j_2} (U\psi_{j_2}U^{*}) U
%\varepsilon_{t-j_2-k} \varepsilon_{t-j_1}^{*}U^{*} (U\psi_{j_1}U^{*})
%\bigg) U \nonumber\\
%&\stackrel{\operatorname{same LSD}}{=}& \sum_{t,j_1,j_2} \Lambda_{j_1}
%\varepsilon_{t-j_1} \varepsilon_{t-j_2-k}^{*}\Lambda_{j_2} +
%\sum_{t,j_1,j_2} \Lambda_{j_2} \varepsilon_{t-j_2-k}
%\varepsilon_{t-j_1}^{*}\Lambda_{j_1}.
%\end{eqnarray}
%% Note that, as $UU^{*}=I$ and for any two square matrices
%$A$ and $B$, $AB$ and $BA$ have identical LSD, matrices in second and
%third line also have same LSD (after scaling by $n$).
%
%Thus it follows that as far as the LSD of sample autocovariance
%matrices are concerned, this is equivalent to the model
%%Now, by (\ref{eqn: diagassump}),
%%we may as well assume that $\{\psi_{j}\}$ are diagonal matrices.
%Under this assumption, $\{\psi_{j,p}=\operatorname{diag}(
%\psi_{j,(1,1)},
%\psi_{j,(2,2)},\ldots,\psi_{j,(p,p)})\}$
%%lead us to the situation where %model without spatial dependence
%(dependence among components). It is because
%and
%
%e1.4 #&#
\begin{equation}
\label{eqn: situation} X_{t,i}\ \mbox{($i$th component of $X_t$)} =
\sum_{j=0}^{\infty} \psi _{j,(i,i)}
\varepsilon_{t-j,i}\qquad\forall i \geq1,
\end{equation}
where $\Lambda_j=\operatorname{diag}(\psi_{j,(1,1)},\psi
_{j,(2,2)},\ldots,\psi
_{j,(p,p)})$ for every $j$.
Hence, this model does not exhibit spatial dependence or dependence
among the components.
%But in practice, we may have both spatial and temporal dependence.

%As we discussed earlier, it is difficult to deal with the LSD of the
%non-symmetric autocovariance matrices $\{\hat{\Gamma}_{k}\}$.
Our approach differs from the existing approaches in many ways. First,
we do away with the Hermitian and simultaneously diagonalizable
condition and replace it with a more natural and much weaker joint
convergence assumption [assumption (A3) in Section~\ref{sec: main}].
Second, all the existing works concentrate on $\hat{\Gamma}_{u} +
\hat
{\Gamma}_{u}^{*}$.
If we wish to study the \textit{singular values} of $\hat{\Gamma}_{u}$,
we need to consider the symmetric product $\hat{\Gamma}_{u}\hat
{\Gamma
}_{u}^{*}$.
This gives rise to a completely different LSD problem. Indeed, one may
consider more general symmetrizations that involve several $\hat
{\Gamma}_{u}$.
As we may recall, in the one-dimensional case, all tests for white
noise are based on quadratic functions of autocovariances. See, for
example, \citet{HL2003}, \citet{S2011} and \citet{XW2014}. The analogous
objects in our model are quadratic polynomials in autocovariances.
Thus, we are naturally led to the consideration of matrix polynomials
of autocovariances. While it is conceivable that the Stieltjes
transform method can be potentially used to tackle these cases, it
seems to be rather cumbersome and needlessly lengthy to do so and shall
at best be a case by case study. We provide a unified method to study
the LSD of symmetric polynomials of the autocovariance matrices. We do
not use Stieltjes transforms at all except to cross-check our results
with the existing results, all of which follow as special cases.
%$\hat{\Gamma}_{k} + \hat{\Gamma}_{k}^{*}$ cannot provide any
%information about the singular values of $\{\hat{\Gamma}_{k}\}$.
%Our general result for the LSD will be applicable to these matrices
%under the same set of assumptions on the coefficient matrices.

%Therefore, in high-dimensional set up, to construct a test for white
%noise and to study the asymptotic distribution of the test statistic,
%one may analogously need to deal with expected average trace of any
%polynomial in %$\{\hat{\Gamma}_{k}\}$. In literature, there is no such
%result regarding the expected average trace of polynomials in $\{\hat{
%\Gamma}_{k}\}$ except a particular one $\{\hat{\Gamma}_{k} + \hat{
%\Gamma}_{k}^{*}\}$.

To obtain the LSD,
% of any symmetric polynomial in $\{\hat{\Gamma}_{u}\}$,
we use the method of moments. The $h$th order moment of the ESD of an
$n \times n$ real symmetric matrix $R_n$ equals
%\begin{equation} \label{eqn: ESDmoments}
$\beta_{h}(R_n): = %\frac{1}{n} \sum_{i=1}^n \lambda_i^h =
\frac{1}{n} \Tr( R_n^h)$.
%\end{equation}
Consider the following conditions:
\begin{longlist}[(M4)]
\item[(M1)] For every $h \ge1, E(\beta_h(R_n)) \rightarrow\beta
_h$, %
%\vskip5pt
%(M2) $\Var(\beta_h (R_n)) \rightarrow0,  \forall h \ge
%1.$

\item[(M4)] $\sum_{n=1}^{\infty}E(\beta_{h}(R_n)- E(\beta
_{h}(R_n)))^4 < \infty, \forall h \geq1$, and

\item[(C)] The sequence $\{ \beta_h\}$ satisfies Carleman's condition,
%\begin{equation}\label{eqn: carle}%\textit{Carleman's Condition:}
%
$ \sum_{h=1}^\infty
\beta_{2h}^{-{1}/{(2h)}} = \infty$.
%\end{equation}
%\begin{lemma} \label{lem:LSDmoment}
\end{longlist}

If (M1), (M4) and (C) hold, then ESD of $R_n$ converges
almost surely to the
distribution $F$ determined uniquely by the moments $\{\beta_h\}$.
%%Further, the convergence is almost sure if $(M4)$ holds.
%\end{lemma}
(M1) is the most crucial condition in this method as it identifies the
moments of the LSD.

In Theorem~\ref{thmm: lsd}, we claim the existence of the LSD of any
symmetric polynomial in $\{\hat{\Gamma}_{u}\}$ in model (\ref
{eqn:model}) and describe the limit in terms of a polynomial of some
free variables. To establish (M1), we use tools from
non-commutative free probability theory (the next section and
Sections \ref{sec: main}, \ref{subsec: ncpconvergence}, \ref{subsec:
lsdproof}
and \ref{subsec: algorithm} contain the necessary background). Free
variables in the non-commutative world are the analogue of independent
random variables in the commutative world. As matrices are
non-commutative objects, appearance of non-commutative probability
spaces is not surprising. The reason for the appearance of free
variables is more subtle (see the discussion at the beginning of
Section~\ref{sec: main}).
In Section~\ref{subsec: example}, we provide simulation results for
specific choices of the model. These simulations support the conclusion
of Theorem~\ref{thmm: lsd}. Based on simulations, we also conjecture
that the LSD exists for the non-Hermitian matrices $\hat{\Gamma}_u$.

It is natural to anticipate that the sample autocovariance matrices
will play an increasingly crucial role in the statistical analysis of
these models. This seems to be at a rudimentary stage currently, but we
anticipate further thrust as the limiting structure of these matrices
is uncovered.
\citet{L2013} estimated the spectrum of the coefficient matrices by
minimizing some distance between Stieltjes transformations of the ESD
and the LSD of $\{\hat{\Gamma}_{u}+\hat{\Gamma}_{u}^{*}\}$ in some
appropriately chosen space of distribution functions.
In Sections \ref{subsubsec: maarorder} and \ref{subsec: arorder}, we
use the LSD results to provide a graphical method to determine the
order of a moving average and an autoregressive process.
%The behaviour of the LSD depends on the order of the process (see
%Remark~\ref{rem: diagnosis}). As discussed before, in Section~\ref{subsubsec: maarorder}, we use this property of sample
%autocovariance matrices to suggest a graphical method for determining
%the order of the moving average and autoregressive processes from a
%data set.
Following a suggestion by one of the referees, in Section~\ref
{subsubsec: trace}, we discuss the asymptotic distribution of the trace
of sample autocovariance matrices. As a by-product of the calculations
used in the derivation of the LSD results, we conclude that these
traces have asymptotic normal distributions. This can be used to test
simple null and alternative hypotheses for model (\ref{eqn:model}).

Section~\ref{sec: proof} contains the outline of the proofs. Further
details of the proofs are available in the supplementary file \citet
{BB2014freesup}.
%This suggests potential statistical applications of the LSD of the
%symmetric polynomials in sample %autocovariance matrices.
%As mentioned before, this suggests a diagnostic method to determine
%the order of the process (see Section~\ref{sec: simulation} for
%details).

%s2 #&#
\section{Some notions from free probability}\label{sec: free}
%In this section, we briefly discuss the concept of free probability
%and provide a summary of tools which are useful for our purpose.

%As discussed in Section~\ref{sec: Intro}, Free variables in the
%non-commutative world, is the analogue of independent random variables
%in commutative world. Note that c
To aid the reader, we first highlight the essential notions from free
probability that are needed to understand our main theorem. Further
concepts and facts, as needed later, are discussed at the beginning of
Section~\ref{sec: main} and in Sections \ref{subsec: ncpconvergence},
\ref{subsec: lsdproof} and \ref{subsec: algorithm}.
An excellent reference for all the details is \citet{NS2006}.

Commutative random variables are attached to a probability space
$(\mathcal{S},E)$, which consists of a $\sigma$-field $\mathcal{S}$ and
an expectation operator $E$. Similarly, non-commutative variables are
attached to some non-commutative $*$-probability space (NCP) $(\mathcal
{A},\varphi)$ consisting of a unital $*$-algebra $\mathcal{A}$ and a
unital linear functional (called a state) $\varphi: \mathcal{A} \to
\mathbb{C}$, $\varphi(1_{\mathcal{A}}) =1$. Thus, $\varphi$ is the
analogue of the expectation operator.
%a variable $a$, $\{\varphi(\Pi(a, a^{*}):  \Pi\
%\operatorname{polynomials}
%\}$ (all moments) is referred as the %distribution of $a$.
The elements of $\mathcal{A}$ are called (non-commutative random) variables.
The canonical example of NCP that we will need is $\mathcal{M}_d$, the
space of all $d\times d$ matrices with the state $\varphi$ as the
average trace. If the matrix has random entries, we modify $\varphi$ by
taking its usual expectation.

In the commutative case, random variables (say with bounded support)
are independent if and only if all joint moments obey the product rule.
It is well known that the cumulants and moments are related via the M\"
{o}bius transformation on the partially ordered set (POSET) of \textit
{all partitions}. Using this, it can be shown that independence is also
equivalent to the vanishing of all mixed cumulants.

For a set of non-commutative variables $\{a_i\}_{i \geq1}$, the set of
all joint moments is defined as $\{\varphi(\Pi(a_i, a_{i}^{*}: i\geq
1)): \Pi \mbox{ polynomials}\}$ and is known as the distribution of
$\{
a_i\}$.
Here, we have the notion of joint cumulants, called \textit{free
cumulants}. These can be uniquely obtained from the above moments and
vice versa via a different M\"{o}bius transformation and its inverse on
the POSET of all \textit{non-crossing partitions}.
Non-commutative variables are said to be \textit{free} (freely
independent) if and only if all their mixed free cumulants vanish.
%Equivalently, let $\varepsilon_{i} = 1$ or $*$ for all $1 \leq i \leq
%r$.
%Variables $a_1,a_2,\ldots, a_r \in\mathcal{A}$ are free if and only
%if their joint free cumulants satisfy
%\begin{equation}
%k_{n}(a_{i_1}^{\varepsilon_{i_1}}, a_{i_2}^{\varepsilon_{i_2}},\ldots,
%a_{i_u}^{\varepsilon_{i_u}}) = 0, \mbox{if $i_j \neq i_{j+1}$ for some
%$1 \leq j \leq u-1, \forall u \geq2$}.
%\end{equation}

A consequence of freeness is that all joint moments of free
variables are computable in terms of the moments of the individual
variables. Of course, the algorithm for computing moments under
freeness is different from (and more complicated than) the product rule
under usual independence.
The notion of freness of variables extends to freeness of sub-algebras
in the natural way. Now consider NCPs $\{(\mathcal{A}_{u},\varphi
_{u})\}
_{1 \leq u \leq r}$. Then, analogous to the product space in the
commutative case, we can have $(\mathcal{A},\varphi)$, the \textit{free
product} of $\{(\mathcal{A}_{u},\varphi_{u})\}$ so that the restriction
of $\varphi$ on $\mathcal{A}_{u}$ is $\varphi_{u}$ and $\mathcal{A}_u$
are free sub-algebras of $\mathcal{A}$.

While matrices are seldom free, there is a large class of matrices that
are free in an asymptotic sense (which is made precise in Section~\ref
{subsec: ncpconvergence}) as the dimension increases. For example, if
$W_1$ and $W_2$ are $n \times n$ independent symmetric matrices with
all entries i.i.d. whose all moments are finite, then they are
asymptotically free. Using such asymptotic freeness, we shall be able
to compute the limits of required traces by using tools from free
probability. This will help us to establish the (M1) condition and in
the bargain also provide us with expressions for the limits in terms of
free variables.
%s3 #&#
\section{Main result} \label{sec: main}

%This section provides the joint convergence of sample
%autocovariance matrices (see Theorem~\ref{thmm: lsd}). As a corollary,
%this implies
%In this section, we establish the existence of the LSD of any
%symmetric polynomials of $\{\hat\Gamma_{i}, \hat\Gamma_{i}^{*}: i \geq
%0\}$, which includes in particular $\hat{\Gamma}_{0}$, $\hat{
%\Gamma}_{i}+\hat{\Gamma}_{i}^{*}$, $\hat{\Gamma}_{i}\hat{
%\Gamma}_{i}^{*}$, $\hat{\Gamma}_{i}\hat{\Gamma}_{i}^{*} + \hat{
%\Gamma}_{j}\hat{\Gamma}_{j}^{*}$ for all $i,j \geq1$. In particular,
%this also implies the results of \citet{JWBN2014} and \citet{LAP2013}.
Consider the following assumptions on the driving process $\{
\varepsilon
_t\}$ and the coefficient matrices $\{\psi_j\}$:
\begin{longlist}[(A1)]
\item[(A1)] $\{\varepsilon_{t,j}\}$ are independent with
$E(\varepsilon_{t,j})=0$ and $E|\varepsilon_{t,j}|^{2} = 1, \forall i,j$.

\item[(A2)] $\sup_{t,j} E(|\varepsilon_{t,j}|^{k})<C_k <
\infty, \forall k \geq1$
or,
for some sequence $\eta_n \downarrow0$, $|\varepsilon
_{t,j}|<\eta_{n}\sqrt{n}, \forall i,j$.%Assumption (A1), (A2),
%(A9) and either (A3) or (A4) hold.
\item[(A3)] $\{\psi_j\}$ are compactly supported and for
any polynomial $\Pi$ in $\{\psi_j, \psi_{j}^{*}\}$, $\lim
p^{-1}\operatorname
{Tr}(\Pi)$ exists and is finite.
\end{longlist}

Later we shall relax assumption (A2).

To see how freeness comes into the picture, and hence how it motivates
the statement of our main theorem, let us focus on $\hat{\Gamma}_{0}$ when
\[
X_t = \varepsilon_t + \psi_{1}
\varepsilon_{t-1.}
\]
Let $Z = (\varepsilon_1, \varepsilon_2,\ldots, \varepsilon_n)_{p
\times n}$ be the independent (ID) matrix. For $i \geq0$, let $P_{i}$
be the $n \times n$ matrix whose $i$th upper diagonal is $1$ and $0$
otherwise. Note that $P_0 = I_n$. For $i <0$, let $P_i = P_{-i}^{T}$ be
the transpose of $P_{-i}$.
Note that
\begin{eqnarray*}
\hat{\Gamma}_{0}&=&n^{-1} \sum
_{t=1}^n (\varepsilon_t +
\psi_{1} \varepsilon_{t-1}) ( \varepsilon_t +
\psi_{1} \varepsilon_{t-1})^{*}
\\
&=& n^{-1} \bigl(ZP_{0}Z^{*} +
\psi_{1}ZP_{0}Z^{*}\psi_{1}^{*}
+ \psi _{1}ZP_{1}Z^{*} + ZP_{-1}Z^{*}
\psi_{1}^{*} \bigr) +R_n
\\
&=& \Delta_{0} +R_n \qquad\mbox{(say)}.
\end{eqnarray*}
%
% It is easy to see (details are provided in Section $1$ of
%Supplementary file \citet{BB2014freesup}) that $p^{-1}
%\mbox{Tr}(R_n^{2}) \stackrel{\operatorname{a.s.}}{\rightarrow} 0$. Now
%for any
%two symmetric $p \times p$ matrices $A$ and $B$ with respective ESD
%$F^{A}$ and $F^{B}$, by $A.41$ of \citet{BS2009},
%\begin{eqnarray}
%L^{3}(F^{A},F^{B}) \leq p^{-1} \operatorname{Tr}(A-B)^{2},
%\end{eqnarray}
%where $L$ is the Levy metric between two distribution functions.
%Therefore,
By Lemma $7.1$ of the supplementary file \citet{BB2014freesup}, $\hat
{\Gamma}_{0}$ and $\Delta_{0}$ have identical LSD.
Thus, our primary goal is to show that for all $r \geq1$, $\lim p^{-1}
E \operatorname{Tr}(\Delta_{0}^{r})$ exists. To achieve this, we
first define
an NCP generated by these matrices. However, the matrices $Z$, $\{
I_{p},\psi_1\}$ and $\{P_0, P_{1}, P_{-1}\}$ are all of different
orders. Therefore, we embed these matrices into larger square matrices
of order $(n+p)$. We embed $Z$ into a Wigner\setcounter{footnote}{1}\footnote{A Wigner matrix
is a square symmetric random matrix with independent mean $0$ variance
$1$ entries on and above the diagonal.} matrix $W$ of order $(n+p)$. Thus,
%
%e3.1 #&#
\begin{equation}
\label{eqn: wembed} W= \pmatrix{ W^{(1)}& Z \cr Z^{*} & W^{(2)}},
\end{equation}
where $W^{(1)}$ and $W^{(2)}$ are two independent Wigner
matrices of order $p$ and $n$ respectively and also independent of $Z$
and whose entries satisfy assumption (A2). For any matrices $B$ and
$D$ of order $p$ and $n$, respectively, let $\bar{B}$ and $\underbar
{D}$ of order $(n+p)$ be the matrices
%
%e3.2 #&#
\begin{equation}
\label{eqn: upperbar}\bar{B} = \pmatrix{ B& 0 \cr 0 & 0}, \qquad\underbar{D} = \pmatrix{ 0 &0 \cr 0 &D}.
%\nonumber
\end{equation}
Note that for any integer $r$, if the right-hand side limits
below exist, then
%
%e3.3 #&#
%e3.4 #&#
%e3.5 #&#
\begin{eqnarray}
\label{eq:upperbar}\lim(n+p)^{-1} \operatorname {Tr} \bigl(
\bar{B}^r \bigr) &=& y(1+y)^{-1} \lim p^{-1}
\operatorname{Tr} \bigl(B^r \bigr),
\\
\label{eq:lowerbar}\lim(n+p)^{-1} \operatorname{Tr} \bigl(
\underbar{D}^r \bigr) &=& (1+y)^{-1} \lim n^{-1}
\operatorname{Tr} \bigl(D^r \bigr) \qquad\mbox{and}
\\
\label{eq:deltabar} \lim p^{-1} \operatorname{Tr} \bigl(
\Delta_{0}^{r} \bigr) &=& y^{-1}(1+y)
\lim(n+p)^{-1} \operatorname{Tr} \bigl(\bar{\Delta}_{0}^{r}
\bigr).
\end{eqnarray}
%
%\begin{equation}\label{eq:lowerbar}\lim(n+p)^{-1} \operatorname{Tr}(
%\underbar{D}^r) = (1+y)^{-1} \lim n^{-1} \operatorname{Tr}(D^r).
%\end{equation}
On the other hand,
\[
n\bar{\Delta}_{0} = \bar{I}_{p} W \underbar{P}_{0}
W \bar{I}_{p} + \bar {\psi}_{1} W \underbar{P}_{0}
W \bar{\psi}_{1}^{*} + \bar{\psi }_{1} W
\underbar{P}_{1} W \bar{I}_{p} + \bar{I}_{p} W
\underbar{P}_{-1} W \bar {\psi}_{1}^{*}.
\]
%
%and provided the right side limit below exists,
%\begin{equation}\label{eq:deltabar} \lim p^{-1} \operatorname{Tr}(
%\Delta_{0}^{r})=y^{-1}(1+y)\lim(n+p)^{-1} \operatorname{Tr}(\bar{
%\Delta}_{0}^{r}).\end{equation}
Thus, $\bar{\Delta}_0^{r}$ involves polynomials in these matrices. So
it is a question of computing the limiting trace of such polynomials.
Now observe that for any monomial $m$:
\begin{longlist}[(1)]
\item[(1)] $\lim p^{-1} \operatorname{Tr}(m({P}_{0}, {P}_{1},{P}_{-1}))$
exists and can be computed easily.
% = (1+y)^{-1} \lim$ $p^{-1}\operatorname{Tr}(m(\underbar{P}_{0},
%\underbar{P}_{1},\underbar{P}_{-1}))$.

\item[(2)] under assumption (A3), $\lim n^{-1} \operatorname{Tr}(m(\bar
{I}_{p}, \bar{\psi}_{1},\bar{\psi}_{1}^{*}))$ exists.

Moreover, from random matrix theory it is well known that

\item[(3)] if (A1) and (A2) hold then $\lim E \operatorname
{Tr}((n+p)^{-1/2}W)^r = E(s^{r})$, where $s$ is a standard semi-circle
variable with moments
%
%e3.6 #&#
\begin{equation}
\label{eqn: smoments} \varphi \bigl(s^{k} \bigr) = %
\cases{
\displaystyle\frac{k!}{({k}/{2})! ({k}/{2}+1)!},&\quad $\mbox{if $k$ is even,}$ \vspace*{2pt}
\cr
0,&\quad $\mbox{if $k$ is odd}$.}
\end{equation}
%
%Some results on NCP and free probability theory guarantee that for any
%monomial $m(\psi_{1}, \psi_{1}^{*})$, $\lim(n+p)^{-1}
%\operatorname{Tr}(m(
%\bar{\psi}_{1},\bar{\psi}_{1}^{*})) < \infty$ is enough for the
%existence of $\lim(n+p)^{-1} \operatorname{Tr}(\bar{\Delta}_{0}^{r})$
%for all
%$r \geq1$ and the limit depends only on the followings.
%\vskip5pt

\item[(4)] Finally, results from free probability guarantee that
\textit{in the limit}, the matrices $(n+p)^{-1/2}W$, $\{\bar
{I}_{p},\bar
{\psi}_1\}$ and $\{\underbar{P}_{0},\underbar{P}_{1},\underbar
{P}_{-1}\}
$ are free variables say $s, \{a_0,a_1\}$ and $\{c_0, c_1, c_{-1}\}$
where $c_1^*=c_{-1}$ in some NCP $(\mathcal{A}, \varphi)$.
\end{longlist}

Thus, using the above conclusions (1), (2) and (3) in conjunction with
equations~(\ref{eq:upperbar}), (\ref{eq:lowerbar}), (\ref{eq:deltabar})
and (\ref{eqn: smoments}), we can conclude that
$\lim p^{-1} \operatorname{Tr}(\Delta_{0}^{r})$ exists and
%can be expressed in terms of appropriate free variables as
%
%e3.7 #&#
\begin{equation}
\label{eqn: 1ma} \lim p^{-1} \operatorname {Tr} \bigl(\Delta
_{0}^{r} \bigr)= y^{-1}(1+y)\varphi \biggl((1+y)
\sum_{j,j^\prime=0,1}a_{j}s c_{j-j^{\prime}}sa_{j}^{*}
\biggr)^{r}.
\end{equation}
%
%??why is there a square in the above formula?. also check next para..
The factor $(1+y)$ within $\varphi$ is the adjustment needed for the
replacement of $Z/\sqrt{n}$ by $W/\sqrt{n+p}$.
The right-hand side of the above equation, involving free variables,
are then the moments of the LSD of $\hat\Gamma_0$.

This is the idea we implement in the general $\operatorname{MA}(q)$ process and for
general symmetric polynomials of the autocovariances. Now we have $q$
coefficient matrices $\{\psi_j\}$ and $\{P_{i}: i=0,\pm1,\pm2,\ldots
\}
$. To describe the limit, consider the NCP $(\mathcal{A},\varphi)$, the
\textit{free product} of the semi-circle variable $s$, $\{a_j\}$ and
$\{
c_i\}$
such that $\varphi(s^k)$ is given by (\ref{eqn: smoments}) and for any
finite monomial $m$, %$m( c_i, c_{-i}, : \ i \geq0)$ and $m(a_j,
%a_j^*, : \ i \geq0 )$,
we have
%
%e3.8 #&#
\begin{equation}
\label{eqn: limitc} \varphi \bigl( m ( c_i, c_{-i} : i \geq0
) \bigr) = (1+y)^{-1}I(J=0), %\begin{cases}
%\frac{1}{1+y} ,\ \ \operatorname{ if $c_i$ and $c_{-i}$ appear equal
%number of
%times %$\# c_i = \# c_{-i},\
%$\forall i \geq1$} \\ %there is same number of $c_i$ and
%%$c_i^*$ for all $i\ge1$} \\
%0,\ \ \operatorname{otherwise}, \end{cases} \nonumber
%J &=& \operatorname{sum of the indices of $\{c_i\}$ appear in $m$.}
%\nonumber
\end{equation}
where $J$ is sum of the subscripts of $\{c_i\}$ which appear in $m$ [the
right-hand side of (\ref{eqn: limitc}) equals $ \lim
n^{-1}\operatorname
{Tr}(m(P_i, P_{-i}: i \geq0))$ %are the limits corresponding to $\{
%\underbar{P}_i,\underbar{P}_i^{*}\}$ and
and can be checked by direct calculation], and
%
%e3.9 #&#
\begin{equation}
\label{eqn: limita} \varphi \bigl( m \bigl(a_j, a_j^* : i
\geq0 \bigr) \bigr) = \frac{y}{1+y} \lim\frac{1}{p} \Tr \bigl( m
\bigl( \psi_j, \psi _j^* : j \geq0 \bigr) \bigr)
\end{equation}
[the right-hand side of (\ref{eqn: limita}) %are the limits
%corresponding to $\{\bar{\psi}_j,\bar{\psi}_j^{*}\}$ which
exists by assumption (A3)].

Let us define for all $u=0,1,2,\ldots,$
%\begin{eqnarray}
%
%e3.10 #&#
\begin{eqnarray}
\label{eqn: gammafree} \gamma_{uq}& =& (1+y)\sum_{j,j^\prime=0}^q
a_j s c_{j^\prime- j+u} s a_{j^\prime}^*,
\nonumber
\\[-8pt]
\\[-8pt]
\nonumber %\nonumber\\
 \gamma_{uq}^* &=& (1+y)\sum_{j,j^\prime=0}^q
a_{j^\prime} s c_{-j^\prime+ j-u} s a_{j}^*. %\gamma_{i\infty} &=& \sum_{j=0}^\infty\sum_{j^\prime=0}^\infty a_j s
%c_{j^\prime- j+i} s a_{j^\prime}^*, \gamma_{i\infty}^*
%= \sum_{j=0}^\infty\sum_{j^\prime=0}^\infty a_{j^\prime} s
%c_{j^\prime- j+i}^* s a_{j}^*.
\end{eqnarray}
%
%\end{eqnarray}
% Next, for all $q \geq0$, consider the NCPs $(
%\mathcal{B}_{q},\varphi)$, %and $(\mathcal{B}_{\infty},\varphi)$,
%where
%\begin{eqnarray} \label{eqn: bq}
%\mathcal{B}_{q} &=& \operatorname{Span}\{(1+y)\gamma_{iq}, (1+y)
%\gamma_{iq}^{*}: i\geq0\} \subset\mathcal{A}.
%%\mathcal{B}_{\infty} &=& \operatorname{Span}\{\gamma_{i\infty},
%\gamma_{i
%\infty}^{*}: i\geq0\} \subset\mathcal{A}.
%\end{eqnarray}
Then we have the following theorem, the proof of which is given in
Section~\ref{subsec: lsdproof}.

%th3.1 #&#
\begin{theorem} \label{thmm: lsd}
Suppose $X_t \sim\operatorname{MA}(q)$, $q < \infty$ and \textup{(A1)}, \textup{(A2)},
\textup{(A3)} and $pn^{-1} \to y \in(0,\infty)$ hold.
Then the LSD of any symmetric polynomial $\Pi(\hat{\Gamma}_{u},\hat
{\Gamma}_{u}^{*}: u \geq0)$
%of $\{\hat{\Gamma}_{u}, \hat{\Gamma}_{u}^{*}\}$
exists and the limit moments are given by
%can be obtained from the relation%the LSD has same distribution as $
%\Pi(\gamma_{iq},\gamma_{i}^{*}: i \geq0)$.
%
%e3.11 #&#
\begin{equation}\qquad
\label{eqn: lsdmoments} \lim p^{-1}E \operatorname{Tr} \bigl(\Pi \bigl(\hat{
\Gamma}_{u},\hat {\Gamma}_{u}^{*}: i \geq0 \bigr)
\bigr) = y^{-1}(1+y)\varphi \bigl(\Pi \bigl(\gamma_{uq},
\gamma _{uq}^{*}: i \geq0 \bigr) \bigr).
\end{equation}
\end{theorem}

%\textbf{????why are yu making this remark? motivation??}
%\begin{remark} \label{rem: jntlsd}
%Let $\Pi_{1}$ and $\Pi_2$ (symmetric) be two polynomials. Then LSD of $
%\Pi_{1}(\psi_{j},\psi_{j}^{*}: j \geq0) \Pi_{2}(\hat{\Gamma}_{u},\hat{
%\Gamma}_{u}^{*}: u \geq0) \Pi_{1}^{*}(\psi_{j},\psi_{j}^{*}: j \geq
%0) $ exists and the limit moments can be obtained from the relation
%\begin{eqnarray} \label{eqn: jointlsd}
% && \lim p^{-1}E \operatorname{Tr}\left(\Pi_{1}(\psi_{j},\psi_{j}^{*}:
%j \geq
%0)\Pi_{2}(\hat{\Gamma}_{u},\hat{\Gamma}_{u}^{*}: i \geq0) \Pi_{1}^{*}(
%\psi_{j},\psi_{j}^{*}: j \geq0)\right) \nonumber\\
% &=& y^{-1}(1+y)\varphi\left(\Pi_{1}(\bar{a}_{j},\bar{a}_{j}^{*}: j
%\geq0)\Pi_{2}(\gamma_{uq},\Gamma_{uq}^{*}: i \geq0) \Pi_{1}^{*}(
%\bar{a}_{j},\bar{a}_{j}^{*}: j \geq0)\right). \nonumber
%\end{eqnarray}
%\end{remark}
For particular symmetric polynomials, the LSD exist under
relaxed moment assumptions. In the next remark, we consider the LSD of
$\{\hat{\Gamma}_{u}+\hat{\Gamma}_{u}^{*}\}$ and $\{\hat{\Gamma
}_{u}\hat
{\Gamma}_{u}^{*}\}$. Its proof, given in Section $5$ of the
supplementary file \citet{BB2014freesup},
%Proof of (b), given in Section $3(b)$ of \citet{BB2014a},
is based on the same truncation arguments as in \citet{JWBN2014} after
some necessary modifications.

%re3.1 #&#
\begin{remark} \label{remark: truncation} Suppose $X_t \sim
\operatorname
{MA}(q)$, $ q < \infty$, and (A1), (A3) and $pn^{-1} \to y \in(0,
\infty)$ hold. Then the following hold true:
\begin{enumerate}[(a)]
\item[(a)] For each $0 \leq u< \infty$, LSD of $\hat{\Gamma}_{u} +
\hat{\Gamma}_{u}^{*}$ exists if for some $\delta\in(0,2]$,
\begin{enumerate}[(A4)]
\item[(A4)] $\sup_{t,j}E(|\varepsilon_{t,j}|^{2+\delta})
< M <\infty$, and

\item[(A5)] for any $\eta> 0$, $\frac{1}{\eta^{2+\delta
}np}\sum_{j=1}^{p}\sum_{t=1}^{n} E(|\varepsilon_{t,j}|^{2+\delta
}I(|\varepsilon_{t,j}|>\eta n^{{1}/{(2+\delta)}})) \rightarrow0$.
%%???range??(A6) and (A7),
\end{enumerate}
\item[(b)] For each $0 \leq u < \infty$, LSD of $\hat{\Gamma
}_{u}\hat{\Gamma}_{u}^{*}$ exists if (A5) holds and
\begin{enumerate}[(A6)]
\item[(A6)] $\sup_{t,j} E|\varepsilon_{t,j}|^{4} < M <
\infty$. %Moreover, if $X_t \sim\operatorname{MA}(q)$, $q \geq0$,
%process and
%Assumptions (A1), (A2), (A7), (A8), (A9) and (A10) hold,
%then
%\vskip5pt
%
%\end{cor}
\end{enumerate}
\item[(c)] Suppose $X_t \sim\operatorname{MA}(0)$ process and assumptions
(A1), (A3) hold. Then the existence of the LSD of $\hat{\Gamma
}_{0}$, under assumption
\begin{enumerate}[(A7)]
\item[(A7)] For any $\eta> 0$, $\eta^{-2}(np)^{-1} \sum_{j=1}^{p}\sum_{t=1}^{n} E(|\varepsilon_{t,j}|^2I(|\varepsilon
_{t,j}| >
\eta\sqrt{n})) \rightarrow0$,
\end{enumerate}
is well known [see \citet{BS2009}]. \citet{JWBN2014}
established the existence of the LSD of any $\hat{\Gamma}_{u} + \hat
{\Gamma}_{u}^{*}$, $u \geq1$, under assumptions (A4) and (A5).
\end{enumerate}
\end{remark}

%?? To prove (a), we use Lemma~\ref{lem:LSDmoment}. Therefore, one
%needs to verify
%By (\ref{eqn: beforecor1}),
%$(M1)$, $(M4)$ ($(M2)$, if we demand only in probability convergence)
%and \textit{Carleman's condition} (C). $(M1)$ follows from (
%\ref{eqn: beforecor1}). Proofs of $(M2)$, $(M4)$ and
%\textit{Carleman's condition} (C) are given in Section $3(a)$ of
%\citet{BB2014a}.

The next remark in particular says that the main result of
\citet{LAP2013} follows from Theorem~\ref{thmm: lsd} and Remark~\ref
{remark: truncation}.

%re3.2 #&#
\begin{remark} \label{rem: liufinal}
(a) Under (A1), (A3), (A4) and (A5), Theorem~\ref
{thmm: lsd} along with Remark~\ref{remark: truncation} provides moments
of the LSD of $\frac{1}{2}(\hat\Gamma_u+\hat\Gamma_u^*)$. These moments
can be used to get the Stieltjes transform $m(z)$ of this LSD as % for
%the LSD of $\hat\Gamma_u+\hat\Gamma_u^*$ as
%
%e3.12 #&#
%e3.13 #&#
%e3.14 #&#
%e3.15 #&#
\begin{eqnarray}
\label{eqn: mmr} m(z) & = &y^{-1}(1+y) \varphi \bigl( \bigl(B(\lambda
,z)-z \bigr)^{-1} \bigr) \qquad\mbox{where}
\\
\label{eqn: kkr} K(z,\theta) & = &y^{-1}(1+y) \varphi \bigl(h(\lambda,
\theta ) \bigl(B(\lambda,z)-z \bigr)^{-1} \bigr),
\\
\label{eqn: hhr}h(\lambda,\theta) &=& \Biggl(\sum_{j=0}^{\infty
}e^{ij\theta
}a_{j}
\Biggr) \Biggl(\sum_{j=0}^{\infty}e^{-ij\theta}a_{j}^{*}
\Biggr), \qquad \lambda= \bigl\{a_j, a_j^{*}: j \geq0
\bigr\},
\\
\label{eqn: bbr} B(\lambda, z) &=& E_{\theta} \bigl( \cos(u\theta) h(
\lambda,\theta) \bigl(1+ y \cos(u \theta) K(z,\theta) \bigr)^{-1}
\bigr),
\end{eqnarray}
and $\theta$ is a $U(0, 2\pi)$ random variable which is commutative
with $\{a_j, a_j^{*}\}$. Details of the arguments is based on a
recursion formula for moments and is given in Section~\ref{subsection:
liufinal}.

%The main Theorem of \citet{LAP2013} provides the LSD of $\hat\Gamma_u+
%\hat\Gamma_u^*$ in terms of its Stieltjes transform.
(b) As discussed in Section~\ref{sec: Intro}, \citet
{LAP2013} proved the existence of the LSD of $\frac{1}{2}(\hat{\Gamma
}_{u} + \hat{\Gamma}_{u}^{*})$ for the model (\ref{eqn:model}). Their
most crucial assumption was the following. % Consider the following
%assumption. % They assume (B) as follows.

(B) $\{\psi_j\}$ are Hermitian and simultaneously
diagonalizable, norm bounded matrices. There are continuous functions
$f_{j}: \mathbb{R} \to\mathbb{R}$ and a unitary matrix $U$ of order
$p$ such that $U\psi_j U^{*}= \operatorname{diag}(f_{j}(\alpha
_1),f_{j}(\alpha
_2),\ldots, f_{j}(\alpha_p))$. ESD of $\{\alpha_1, \alpha_2,\ldots,
\alpha_{p}\}$ converges weakly to a compactly supported probability
distribution $F_{a}$.

Note that assumption (B) implies assumption (A3). The
main theorem of \citet{LAP2013}, under (B), provides the LSD of
$\frac
{1}{2}(\hat\Gamma_u+\hat\Gamma_u^*)$ with its Stieltjes transform
satisfying %Under Assumption (B), it is easy to see that the
%Stieltjes transform equations (\ref{eqn: mmr})-(\ref{eqn: bbr}) reduce
%to %Under (B), \citet{LAP2013} established the Stieltjes transform
%$m(z)$ for the LSD of $\hat\Gamma_u+\hat\Gamma_u^*$ as
%
%e3.16 #&#
%e3.17 #&#
%e3.18 #&#
\begin{eqnarray}\qquad\qquad
\label{eqn: mr} m(z) %\nonumber\\
&=& \int \biggl( E_{\theta^\prime} \biggl(
\frac{\cos(u\theta^\prime)
h_{1}(\alpha,\theta^\prime) }{1+ y \cos(u \theta^\prime)
K(z,\theta
^\prime)} \biggr)-z \biggr)^{-1} \,dF_{a}(\alpha)\qquad
\mbox{where}%
%\nonumber
\\
\label{eqn: kr} K(z,\theta) %\nonumber\\
&=& \int \biggl( E_{\theta^\prime} \biggl(
\frac{\cos(u\theta^\prime)
h_{1}(\alpha,\theta^\prime) }{1+ y \cos(u \theta^\prime)
K(z,\theta
^\prime)} \biggr)-z \biggr)^{-1} h_{1}(\alpha,
\theta) \,dF_{a} (\alpha ), %
%\nonumber\\
\\
\label{eqn: hr}
h_{1}(\alpha, \theta) &=& \Biggl|\sum_{j=0}^{q}
e^{ij\theta}f_{j}(\alpha)\Biggr|^{2}. % \mbox{and} \alpha\sim
%F_a.
\end{eqnarray}
%
% Give the assumptions and the results ???
It can be shown that under assumption (B), the Stieltjes transform
equations (\ref{eqn: mmr})--(\ref{eqn: bbr}) reduce to equations
(\ref
{eqn: mr})--(\ref{eqn: hr}). %The main Theorem of \citet{LAP2013}, under
%(B), provides the LSD of $\hat\Gamma_u+\hat\Gamma_u^*$ with its
%Stieltjes transform satisfying
% Now since, Under (B), the Stieltjes transform
Thus, Theorem~\ref{thmm: lsd} in conjunction with Remark~\ref{remark:
truncation} implies the main theorem of \citet{LAP2013}. % for the LSD
%of $\hat\Gamma_u+\hat\Gamma_u^*$.
%Under the assumptions of \citet{LAP2013}, the Stieltjes transformation
%obtained there can be easily obtained %using (\ref{eqn: moment01}) and
%the power series expansion of Stieltjes transformation. %given in %(
%\ref{eqn: expansionstieltjes}).
%The Stieltjes transformation in \citet{LAP2013} can easily be derived
%from (\ref{eqn: moment01}).
%However, (\ref{eqn: moment01}) remains true more generally if %the
%coefficient matrices are symmetric \textbf{????} and
%(A1), (A3), (A4), (A5) ???How can you get all moments with
%(A5)? hold.
\end{remark}

So far, we have assumed $q< \infty$. With some additional
assumptions, the results continue to hold for $q=\infty$. The proof of
Corollary~\ref{rem: mainfty} is based on truncation arguments and is
given in Section $6$ of the supplementary file \citet{BB2014freesup}.

%co3.1 #&#
\begin{cor} \label{rem: mainfty}
Theorem~\ref{thmm: lsd}, Remark~\ref{remark: truncation}\textup{(a), (b)}
and Corollary~\ref{rem: final} hold for $\operatorname{MA}(\infty)$ process also,
after replacing $q$ by $\infty$ provided

\textup{(A8)} $\sum_{j=0}^{\infty} \sup_{p} \Vert\psi_j\Vert <
\infty$, where for all $j \geq0$, $\Vert\psi_j\Vert$ denotes maximum
absolute eigenvalue of $\psi_j$. %$ \sum_{j=1}^\infty
%\left\{ \varphi( a_j^*
%a_j a_j^* a_j)\right\}^{\frac{1}{4}} < \infty$ %and for every $
%\varepsilon>0$, there exists a $p_{0} \ge1$ such that
%for every $p \ge p_0$, %there is $q_{1} \ge1$ such that $q \ge q_1$
%%implies
%$\sum_{j=1}^\infty\left\{ p^{-1} \Tr( \psi_j^*
%\psi_j \psi_j^* \psi_j)\right\}^{\frac{1}{4}} < \infty.$
%and there exists $p_0 \geq1$ such that \\ $\lim_{q \to\infty}
%\sup_{p \geq p_0} \sum_{j \geq q} \Vert\psi_j\Vert= 0$.
\end{cor}

%\begin{remark} \label{rem: diagnosis} Theorem~\ref{thmm: lsd} can be
%used for statistical purposes. See
%Sections \ref{sec: application}.
%\ref{subsubsec: maarorder} and \ref{subsec: arorder} for details.
%\end{remark}

%s3.1 #&#
\subsection{Examples}
%
%ex1 #&#
\begin{example} \label{rem: marchenko} Consider the $\operatorname{MA}(0)$ process,
that is, $X_{t} = \varepsilon_{t}\ \forall t$ and suppose assumptions
(A1), (A2), (A3) and $pn^{-1} \to y \in(0,\infty)$ %either (A3) or
%(A4)
hold. Then the following results (a)--(c) follow from Theorem~\ref{thmm:
lsd} and Remark~\ref{rem: liufinal}. %We can relax assumption (A4) and assume
%appropriate weaker conditions on moments (see Remark~\ref{rem:
%truncation}). % and finite $4$-th moment of $\varepsilon_t$ (see
%Section~\ref{subsection: truncation}).

(a) \textit{Mar\v{c}enko--Pastur law}: The LSD of $\hat{\Gamma}_{0}$ is
the Mar\v{c}enko--Pastur law, whose moment sequence is given by [see,
e.g., \citet{MP1967} or \citet{BS2009}]
%
%e3.19 #&#
\begin{equation}
\label{eqn: MPmoments} \beta_{h} = \sum_{k=1}^{h}
\frac{1}{k}\pmatrix{h-1 \cr k-1} \pmatrix{h \cr k}y^{k-1},\qquad  h \geq1.
\end{equation}
%
% Moreover, the result continues to hold if we assume (A5) instead
%of (A3) or (A4).

(b) %\label{rem: bassel}
\textit{Free Bessel law}:
%Recall that the $h$-th order raw moments of Bessel$(2, y^{-1})$ is
%given by (\ref{eqn: besselmoments}). Under the same model as in Remark~\ref{rem: marchenko}, it turns out that
The LSD of $ (\frac{n}{p} )^2\hat{\Gamma}_{u}\hat{\Gamma
}_{u}^{*}$, $u \geq1$ is the free Bessel$(2,y^{-1})$ law,
characterized by the moment sequence,
%
%e3.20 #&#
\begin{equation}
\label{eqn: besselmoments} \beta_{h} = \sum_{k=1}^{h}
\frac{1}{k}\pmatrix{h-1 \cr k-1} \pmatrix{2h \cr k-1} y^{-k},\qquad h \geq1.
\end{equation}
%
%For $y=0$, the LSD of $\left(\frac{n}{p}\right)\hat{
%\Gamma}_{i}(\varepsilon)\hat{\Gamma}_{i}^{*}(\varepsilon)$ is the Mar
%\v{c}enko-Pastur law with parameter $1$.
%Moreover, the LSD is identical for all $i \geq1$.
%If we assume (A7) and (A8) instead of (A3) or (A4), then also
%these results holds. %If $p \to\infty$ such that $y=0$, then
%%following the arguments leading to Remark~\ref{rem: bassel},
%the $h$-th order moment of the LSD of $\left(\frac{n}{p}\right)\hat{
%\Gamma}_{i}(\varepsilon)\hat{\Gamma}_{i}^{*}(\varepsilon)$ is $
%\frac{1}{h}{2h \choose h-1}$. Hence, the corresponding LSD is the Mar
%\v{c}enko-Pastur law with parameter $1$.}

(c) %\label{rem: unfoldmp}
%In this remark,
The LSD of $\frac{1}{2}({\hat\Gamma}_u + {\hat\Gamma}^{*}_u)$ are
identical for all $u \geq1$ and their common Stieltjes transformation
$m(z)$ %of the LSD of $\frac{1}{2}({\hat\Gamma}_i + {\hat\Gamma}^{*}_i)$
satisfies the bi-quadratic equation (with one valid solution)
%
%e3.21 #&#
\begin{equation}
\label{eqn: bai} \bigl(1-y^{2}m^{2}(z) \bigr)
\bigl(yzm(z)+y-1 \bigr)^2 = 1.
\end{equation}
This is Theorem $2.1$ of \citet{LAP2013} for the $\operatorname{MA}(0)$ case
and Theorem~$1.1$ of \citet{JWBN2014}. % Also this result holds if we
%assume (A6) and (A7) instead of (A3) or (A4).}
\end{example}

%As discussed in Remark~\ref{remark: truncation}, Corollary
By Remark~\ref{remark: truncation}, Example~\ref{rem: marchenko}(a)
continues to hold if we assume (A7) instead of~(A2).
If we assume (A5) and (A6) instead of (A2), then Example~\ref
{rem: marchenko}(b) continues to hold.
Moreover, Example~\ref{rem: marchenko}(c) holds if we assume (A4)
and (A5) instead of (A2). %by Corollary~\ref{rem: sympoly}(b),
Justification for Example~\ref{rem: marchenko} is given in
Section~\ref{subsection: marchenko}. % under the Assumption (A3) or
%(A4).
%Relaxation of these assumptions follows by Corollary~\ref{rem:
%sympoly} (b).
%$$((1-y)-zym(z))^{2}(1-y^{2}m^{2}(z)) = 1.$$
%Moreover, this LSD is identical for all $i \geq1$.
%\begin{remark}
%LSD of the sample covariance matrix without independence structures in
%columns has been established in \citet{BZ2008}.
%\end{remark}
%%Hence Theorem $(2.1)$ in \citet{LAP2013} is verified by
%Theorem~\ref{thmm:lsd} (c). }

%ex2 #&#
\begin{example} \label{rem: baiprelim}
%\textit{$\textbf{(c)}$ %\label{rem: stieltjes1}
%%Consider the model $X_{t} = \varepsilon_{t}$ where $\varepsilon_{t}
%\sim IID(0,\Sigma_{p})$ as considered in \citet{BZ2008}.
%
Let $X_t = \varepsilon_t$ where $\{\varepsilon_{t,j}\}$'s are all
i.i.d. random variables with mean $0$, variance $1$ and $E|\varepsilon
_{1.1}|^{2+\delta} < \infty$ for some $\delta>0$. Moreover, suppose
$pn^{-1} \to y \in(0,\infty)$ holds. % If the columns of the $p \times
%n$ matrix $Z_{0}=(e_1,e_2,\ldots,e_n)$
%are i.i.d. with dispersion matrix $\Sigma_p$, t
Then, using the same idea as in the proof of Example~\ref{rem:
marchenko}(c), it can be shown that the Stieltjes transform of the
LSD of $\Sigma^{1/2}\hat{\Gamma}_{0}\Sigma^{1/2}$ is given by %(
%\ref),
%
%e3.22 #&#
\begin{equation}
\label{eqn: baistieltjes} m(z) = \int\frac{dF_{\Sigma}(t)}{z-t(1-y-yzm(z))},
\end{equation}
where $\Sigma$ is a symmetric positive definite matrix with compactly
supported LSD $F_{\Sigma}$.
This is Theorem $1.1$ of \citet{S1995}. If $F_{\Sigma} = \delta_{1}$,
this reduces
to the Stieltjes transform of the Mar\v{c}enko--Pastur law. % (see for
%example %\citet{NS2006}, \citet{CD2011} or
%\citet{BS2009}). %(\ref{eqn: mpstieltjes}).
\end{example}

%For a more general result see \citet{BZ2008}.
%Proof of Example~\ref{rem: baiprelim} follows the same idea and is
%omitted. % is given in Section $8$ of the supplementary file. This
%follows using the same idea used in the proof of Example~\ref{rem:
%marchenko}(c). %(see Section~\ref{subsection: marchenko}(c)).
Apart from the expression (\ref{eqn: lsdmoments}) in terms of free
variables, in general, there is no further simplified form of the LSD
of $(\hat{\Gamma}_{u} + \hat{\Gamma}_{u}^{*})$. In the special case
$\psi_{j} = \lambda_{j}I_{p}$, $\lambda_j \in\mathbb{R}$, for all $j
\geq0$, we can describe the LSD in terms of a compound free Poisson
distribution. We need some preparation for this description.
%Next we shall define \textit{compound free Poisson} distribution. This
%is useful to describe the LSD of $\{\frac{1}{2}(\hat{\Gamma}_{i} +
%\hat{\Gamma}_{i}^{*})\}_{i\geq0}$ (in Corollary~\ref{rem: consta}),
%when the coefficient matrices are $\psi_{j} = \lambda_{j}I_{p}$, for
%all $j \geq0$

%de3.1 #&#
\begin{definition} \label{def: freepoissondefn}
A probability measure $\mu$ on $\mathbb{R}$ with free cumulants
\[
k_{n}(\mu) = \lambda m_{n}(\nu) \qquad\forall n \geq1,
\]
for some $\lambda>0$ and some compactly supported probability measure
$\nu$ on $\mathbb{R}$ with moments $\{m_n(\nu)\}$, is called
a \textit{compound free Poisson distribution} with rate $\lambda$ and
jump distribution $\nu$.
\end{definition}

%Let $(\mathcal{A},\varphi)$ be a non-commutative probability space.
%Let $s, a \in\mathcal{A}$ be such that
As an example, suppose $s$ is a semi-circular variable, defined by the
moment sequence (\ref{eqn: smoments}), and $a$ is another variable free
of $s$. Then the free cumulants of $sas$ are given by [see Proposition
$12.18$ in \citet{NS2006}]
%
%e3.23 #&#
\begin{equation}
\label{eqn: cumsas} k_{n}(sas,sas,\ldots,sas) = \varphi
\bigl(a^n \bigr) \qquad\forall n \geq1.
\end{equation}
In particular, if $a$ is self-adjoint with distribution $\nu$, then
$sas$ has the compound free Poisson distribution with rate $\lambda=
1$ and jump measure $\nu$.

%Let $Z = ((\varepsilon_{i,j}))$ be any $p\times n$ ID matrix
%satisfying (A1) and (A2).
Let $A_{p\times p}$ be self-adjoint with compactly supported LSD $a$.
%Here $p=p(n)$ is such that $pn^{-1} \rightarrow y \in(0,\infty)$.
Then it can be shown that, under (A1) and (A2),
the limiting free cumulants of $ZAZ^{*}$ are given by
%
%e3.24 #&#
\begin{equation}
\label{eqn: zaz} \lim_{n}k_{r}
\bigl(ZAZ^{*},ZAZ^{*}, \ldots,ZAZ^{*} \bigr) =
y^{r-1}\varphi \bigl(a^r \bigr)\qquad \forall r\geq1.
\end{equation}
Therefore, asymptotically $ZAZ^{*}$ is a compound free Poisson variable
with rate $y^{-1}$ and jump distribution $ya$.

Now we are ready to state the next example. Justification for Example~\ref{rem: consta} and Remark~\ref{rem: paffel} are given respectively
in Sections \ref{subsection: constaproof} and \ref{subsection: paffel}.

%ex3 #&#
\begin{example} \label{rem: consta}
Let $X_t \sim\operatorname{MA}(q)$ process and suppose assumptions (A1),
(A2) and $pn^{-1} \to y \in(0,\infty)$ hold. Let $\psi_{j} =
\lambda
_{j}I_{p}$, $1 \leq j \leq q$. Then the LSD of $\frac{1}{2}(\hat
{\Gamma
}_{u}+\hat{\Gamma}_{u}^{*})$ is a compound free Poisson
%(see Definition~\ref{def: freepoissondefn})
whose $r$th order free cumulant equals
%
%e3.25 #&#
\begin{equation}
\label{eqn: constacum} k_{ur} = y^{r-1}E_{\theta} \bigl(
\cos(u \theta)\tilde{h}(\lambda,\theta ) \bigr)^r\qquad \forall i \geq0,
% \operatorname{see }.
\end{equation}
where
%
%e3.26 #&#
\begin{eqnarray}
\label{eqn: h} \tilde{h}(\lambda, \theta) &=& \Biggl| \sum_{j=0}^{q}
e^{i j \theta} \lambda _{j}\Biggr|^{2}, \qquad\lambda_0=1,\qquad
\lambda= (\lambda_1, \lambda _2,\ldots,
\lambda_q) \quad\mbox{and}
\nonumber
\\[-8pt]
\\[-8pt]
\nonumber
\theta&\sim& U(0,2\pi).
\end{eqnarray}
\end{example}

%
%f1 #&#
\begin{figure}[b]

\includegraphics{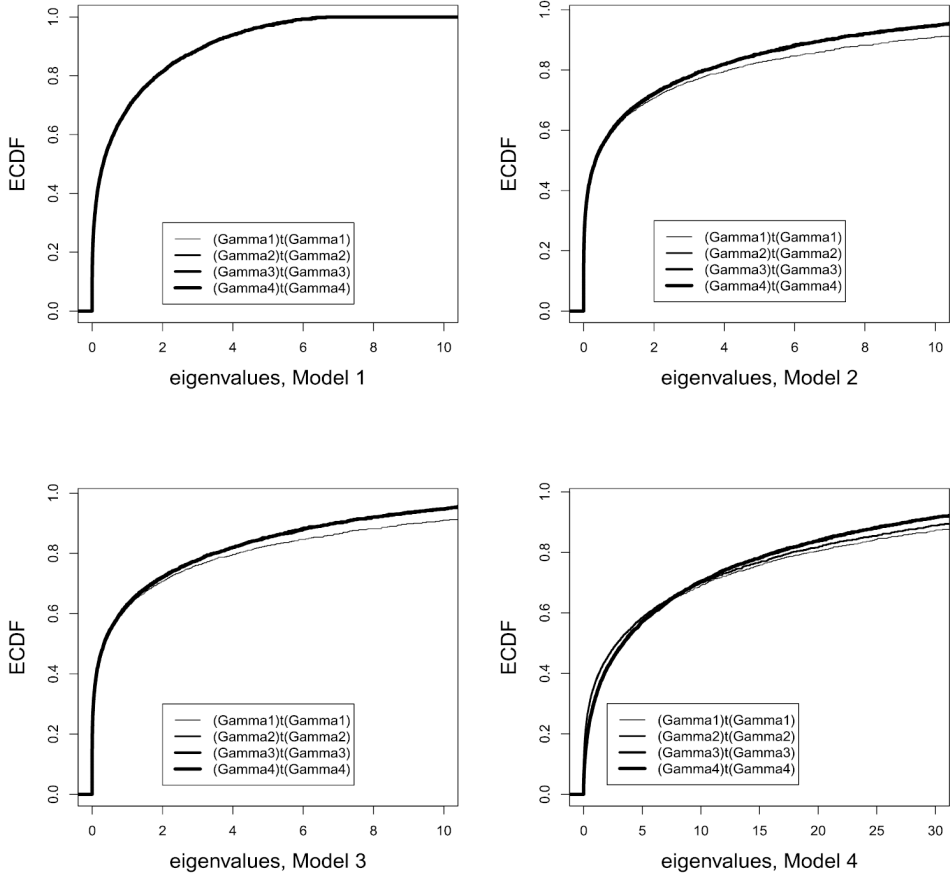}

\caption{ECDF of $\hat{\Gamma}_{u}\hat{\Gamma}_{u}^{*}$, $1 \leq u
\leq4$ for
$n=p=300$.}\label{fig: 1}
\end{figure}

%re3.3 #&#
\begin{remark} \label{rem: paffel}
By Remark~\ref{remark: truncation}(b), (\ref{eqn: constacum})
continues to hold if we assume (A4) and (A5) instead of (A2).
Example~\ref{rem: consta} together with Remark~\ref{remark:
truncation}(b) justifies Theorem $2.1$ in \citet{LAP2013} when $\psi_j
= \lambda_j$, Theorem $1.2$ in \citet{PS2011} and Theorem $1$ in \citet
{Y2012}, though none of them had identified the limit as a
compound free Poisson.
\end{remark}
%
%This of course will be compatible (\ref{eqn: simplifystieltjes}).
% and fourth moment of $\varepsilon_t$ is finite.
% \pagebreak

%
%s4 #&#
\section{Numerical examples and applications} \label{sec: simulation}
%In this section, we first discuss some examples and provide
%simulations to support Theorem~\ref{thmm: lsd} for those examples.
%Then, we discuss some applications of results obtained in Section~\ref{sec: main}. First we see that how to determine order of the
%underlying moving average process from a given data. Under some
%additional restrictions on the parameter space, we also provide a
%method to determine %check white noise for residuals of the fitted
%autoregressive process. This helps to determine
%the order of the autoregressive process. In Section~\ref{subsubsec:
%trace}, as suggested by one of the referees, we discuss the asymptotic
%distribution of traces under model (\ref{eqn:model}). This is not the
%direct application of Theorem~\ref{thmm: lsd}, its corollaries and
%remarks, but an application of Lemma $4.1$ stated and proved in
%Supplementary file.

%s4.1 #&#
\subsection{Numerical examples} \label{subsec: example}
Let $I_p$ and $J_p$ be respectively the identity matrix of order $p$
and the $p \times p$ matrix with all entries $1$ and let $\varepsilon_t
\sim\mathcal{N}_{p}(0, I_{p}), \forall t$. Let $A_{p} = 0.5I_{p}$,
$B_{p} = 0.5(I_{p} + J_{p})$. Let $C_p = ((c_{i,j}))$ and $D_p =
((d_{i,j}))$ be two\vadjust{\goodbreak} $p \times p$ matrices with $c_{i,i} = I(1 \leq i
\leq[p/2])-I([p/2]<i\leq p)$, $d_{i,p+1-i} = 1$ for all $i \geq1$ and
$0$ otherwise. We consider the following models.

\textit{Model} $1$: $X_t = \varepsilon_t$.

\textit{Model} $2$: $X_t = \varepsilon_t + A_{p}\varepsilon_{t-1}$.

\textit{Model} $3$: $X_t = \varepsilon_t + B_{p} \varepsilon_{t-1}$.

\textit{Model} $4$: $X_t = \varepsilon_t + C_{p} \varepsilon
_{t-1} + D_{p} \varepsilon_{t-2}$.

Note that in model $4$, $C_pD_p \neq D_pC_p$, and hence they
are not simultaneously diagonalizable and % (as these two matrices are
%not commutative i.e. $C_pD_p \neq D_pC_p$). As discussed in Section~\ref{sec: Intro} (see (\ref{eqn: ma1example})),
the result of \citet{LAP2013} is not applicable.
For each of these models, we draw one random sample of size
$n$ ($n=300,500$ and $1000$). For each $1 \leq u \leq4$,\vspace*{1pt} we plot the
cumulative distribution function of ESD (ECDF) of $\hat{\Gamma
}_{u}\hat
{\Gamma}_{u}^{*}$ and $\hat{\Gamma}_{u}\hat{\Gamma}_{u}^{*}+ \hat
{\Gamma
}_{u+1}\hat{\Gamma}_{u+1}^{*}$. The graphs for $n=300$ are given in
Figures \ref{fig: 1} and \ref{fig: 2}. Figures $1$ and $2$ in the
supplementary file \citet{BB2014freesup} contain graphs for $n=500$ and
$1000$. These graphs support the following points:

%f2 #&#
\begin{figure}

\includegraphics{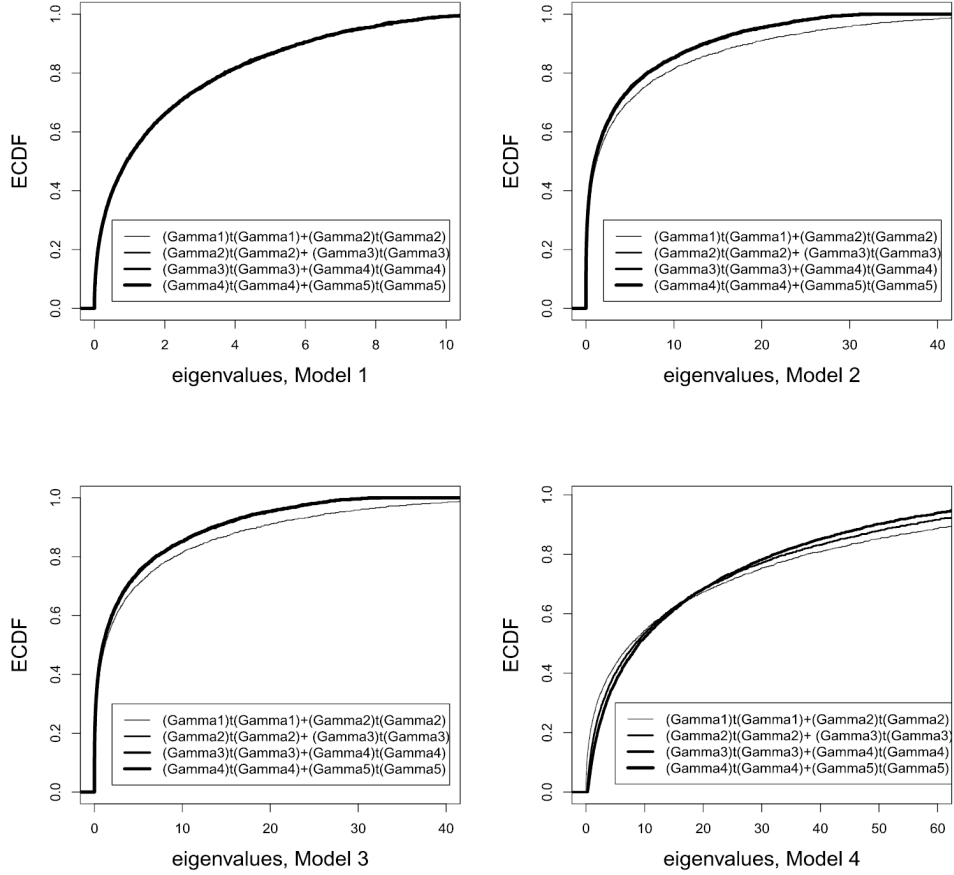}

\caption{ECDF of $\hat{\Gamma}_{u}\hat{\Gamma}_{u}^{*}+ \hat
{\Gamma}_{u+1}\hat
{\Gamma}_{u+1}^{*}$, $1 \leq u \leq4$
for $n=p=300$.} \label{fig: 2}
\end{figure}

% For each of the following model we draw three random samples
%of size $300$, $500$ and $1000$. We consider the dimension $p$ equals
%the sample size $n$. For each of the above models and sample sizes, we
%compute sample autocovariance matrices $\hat{\Gamma}_{1}$, $\hat{
%\Gamma}_{2}$, $\hat{\Gamma}_{3}$ and $\hat{\Gamma}_{4}$ and plot
%cumulative distribution function of ESD of $\hat{\Gamma}_{i}\hat{
%\Gamma}_{i}^{*}$, for $i=1,2,3,4$ in the same graph (see Figure~\ref{fig: 1}). We also do the same for $\hat{\Gamma}_{i}\hat{
%\Gamma}_{i}^{*}+ \hat{\Gamma}_{i+1}\hat{\Gamma}_{i+1}^{*}$,
%$i=1,2,3,4$ (see Figure~\ref{fig: 2}). Simulation results, reported in
%Figures \ref{fig: 1} and \ref{fig: 2}, support the following points.

(a) %Convergence rate is high enough and s
For each of the above models, the ECDF are nearly identical for
$n=300,500$ and $1000$, that is, convergence has already occurred at $n=300$.
For smaller values of $n$, convergence did not occur in our simulation.
Some modification may improve the situation for smaller sample sizes.
Here, we did not investigate any possible modifications.
%Sample size $300$ is enough to apply the asymptotic results discussed
%in Section~\ref{sec: main}. This is because, three columns in Figures
%\ref{fig: 1} and \ref{fig: 2} show identical results, i.e. convergence
%has already occurred at sample size $300$. Sample size less than $300$
%is not showing very good results. Therefore, for small sample size, we
%may need some modifications. We are not going to this investigation in
%this paper.

(b) ECDF of $\hat{\Gamma}_{u}\hat{\Gamma}_{u}^{*}$ (or $\hat
{\Gamma}_{u}\hat{\Gamma}_{u}^{*}+ \hat{\Gamma}_{u+1}\hat{\Gamma
}_{u+1}^{*}$) are almost identical---for all $u > 0$ in model $1$, for
all $u > 1$ in models $2$ and $3$ and for all $u > 2$ in model $4$.
Moreover, ECDFs are different---for $u=1,2$ in both models $2$ and $3$,
and for $u=1,2,3$ in model $4$.
%in both Models $2$ and $3$, ECDFs % of $\hat{\Gamma}_{i}\hat{
%\Gamma}_{i}^{*}$ (or $\hat{\Gamma}_{i}\hat{\Gamma}_{i}^{*}+ \hat{
%\Gamma}_{i+1}\hat{\Gamma}_{i+1}^{*}$)
%are different for $i=1,2$. For $i=1,2,3$, ECDF of $\hat{\Gamma}_{i}
%\hat{\Gamma}_{i}^{*}$ (or $\hat{\Gamma}_{i}\hat{\Gamma}_{i}^{*}+ \hat{
%\Gamma}_{i+1}\hat{\Gamma}_{i+1}^{*}$) are different in Model $4$.
%(b) LSD of $\hat{\Gamma}_{i}\hat{\Gamma}_{i}^{*}$ (or $
%\hat{\Gamma}_{i}\hat{\Gamma}_{i}^{*}+ \hat{\Gamma}_{i+1}\hat{
%\Gamma}_{i+1}^{*}$) are identical for all $i > 0$ in Model $(1)$, for
%all $i > 1$ in Models $(2)$ and $(3)$, for all $i > 2$ in Model $(4)$.
%This is because in Figures \ref{fig: 1} and \ref{fig: 2}, $(i)$ black,
%red, green and blue lines all coincide in Row $(1)$ (Model $(1)$),
%$(ii)$ red, green, blue lines coincide in Row $(2)$ and $(3)$ (Model
%$(2)$ and $(3)$), and $(iii)$ green and blue lines coincide in Row
%$(4)$ (Model $(4)$). Moreover, LSD of $\hat{\Gamma}_{i}\hat{
%\Gamma}_{i}^{*}$ (or $\hat{\Gamma}_{i}\hat{\Gamma}_{i}^{*}+ \hat{
%\Gamma}_{i+1}\hat{\Gamma}_{i+1}^{*}$) are different $(i)$ for $i=1,2$
%in each Model $(2)$, $(3)$ (as black line, in Rows $(2)$ and $(3)$ of
%Figures \ref{fig: 1} and \ref{fig: 2}, is separate from other lines)
%and $(ii)$ for $i=1,2,3$ in Model $(4)$ (as black, red and green
%lines, in Row $(3)$ of Figures \ref{fig: 1} and \ref{fig: 2}, are
%separate from each other).
%This supports Remark~\ref{rem: diagnosis}.

(c) For the $\operatorname{MA}(1)$ process, LSD of $\hat{\Gamma}_{u}\hat
{\Gamma}_{u}^{*}$ (or $\hat{\Gamma}_{u}\hat{\Gamma}_{u}^{*}+ \hat
{\Gamma
}_{u+1}\hat{\Gamma}_{u+1}^{*}$) depends on $\psi_{1}$ only through its
LSD. Since LSD of $A_p$ and $B_p$ are identical (both have mass $1$ at
$0.5$), the ECDF for models $2$ and $3$ are almost identical. % to
%those for Model $3$. %and therefore Row $2$ and $3$ in Figure~\ref{fig: 1} (or Figure~\ref{fig: 2}) show more or less same results.

(d) As noted above, the result of \citet{LAP2013} is not applicable for
model $4$. However, by Theorem~\ref{thmm: lsd}, %in particular, states
%that, under Model $(4)$,
the LSD of any symmetric polynomial in $\{\hat{\Gamma}_{u}, \hat
{\Gamma
}_{u}^{*}\}$ for model $4$ exists and this is supported by row $2$
right panel of Figures \ref{fig: 1} and \ref{fig: 2}.

In Table \ref{table: 1}, we have recorded the mean and variance of the
ESD of $\hat{\Gamma}_{u}\hat{\Gamma}_{u}^{*}$ and $\hat{\Gamma
}_{u}\hat
{\Gamma}_{u}^{*}+\hat{\Gamma}_{u+1}\hat{\Gamma}_{u+1}^{*}$, $1
\leq u
\leq4$, %(or $\hat{\Gamma}_{i}\hat{\Gamma}_{i}^{*}+ \hat{\Gamma}_{i+1}
%\hat{\Gamma}_{i+1}^{*}$)
for model $4$ and $n=p=300$ along with the mean and the variance of
their LSD %of $\hat{\Gamma}_{u}\hat{\Gamma}_{u}^{*}$ and $\hat{
%\Gamma}_{u}\hat{\Gamma}_{u}^{*}+\hat{\Gamma}_{u+1}\hat{
%\Gamma}_{u+1}^{*}$ %(or $\hat{\Gamma}_{i}\hat{\Gamma}_{i}^{*}+ \hat{
%\Gamma}_{i+1}\hat{\Gamma}_{i+1}^{*}$)
using the description of the limits, given in Theorem~\ref{thmm: lsd},
in terms of free variables and limits of coefficient matrices $C_p$ and
$D_p$. The empirical results agree with the theoretical results.

%
%t1 #&#
\begin{table}
\caption{Means and variances for model $4$, $n=p=300$}\label{table: 1}
\begin{tabular*}{\textwidth}{@{\extracolsep{\fill}}ld{2.2}d{2.0}d{2.4}c@{}}
\hline
\textbf{Matrix} & \multicolumn{1}{c}{\textbf{Sample mean}} &
\multicolumn{1}{c}{\textbf{Mean of LSD}} &
\multicolumn{1}{c}{\textbf{Sample variance}} &
\multicolumn{1}{c@{}}{\textbf{Variance of LSD}}
\\
\hline
$\hat{\Gamma}_{1}\hat{\Gamma}_{1}^{*}$ & 10.92 & 11 & 277.8869 &
278 \\
$\hat{\Gamma}_{2}\hat{\Gamma}_{2}^{*}$ & 9.93 & 10 & 214.437 & 215\\
$\hat{\Gamma}_{3}\hat{\Gamma}_{3}^{*}$ & 8.91 & 9 & 143.2908 & 143\\
$\hat{\Gamma}_{4}\hat{\Gamma}_{4}^{*}$ & 8.89 & 9 & 143.2524 & 143\\
$\hat{\Gamma}_{1}\hat{\Gamma}_{1}^{*}+\hat{\Gamma}_{2}\hat
{\Gamma
}_{2}^{*}$ & 20.85 & 21 & 802.6798 & 805\\
$\hat{\Gamma}_{2}\hat{\Gamma}_{2}^{*}+\hat{\Gamma}_{3}\hat
{\Gamma
}_{3}^{*}$ & 18.85 & 19 & 547.4531 & 546\\
$\hat{\Gamma}_{3}\hat{\Gamma}_{3}^{*}+\hat{\Gamma}_{4}\hat
{\Gamma
}_{4}^{*}$ & 17.81 & 18 & 433.1116 & 434\\
$\hat{\Gamma}_{4}\hat{\Gamma}_{4}^{*}+\hat{\Gamma}_{5}\hat
{\Gamma
}_{5}^{*}$ & 17.76 & 18 &433.507 & 434\\
\hline
\end{tabular*}
\end{table}

\begin{figure}

\includegraphics{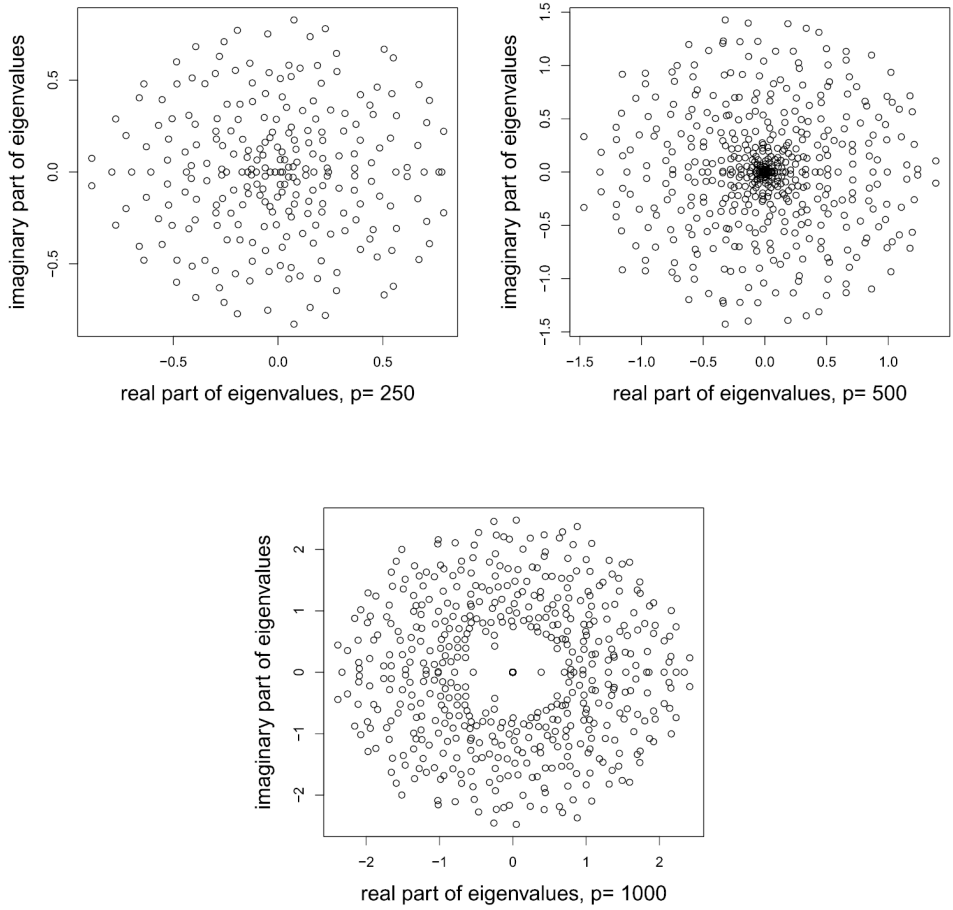}

\caption{ESD of $\hat{\Gamma}_{1}$ for $\operatorname{MA}(0)$ standard
Gaussian process for $n=500$ (Multiple eigenvalues are plotted only
once).}
\label{fig: 5}
\end{figure}

Incidentally, the autocovariance matrices $\{\hat{\Gamma}_{u}\}$
themselves are not symmetric for $u \geq1$ and Theorem~\ref{thmm: lsd}
does not apply. Nevertheless, their ESD should also converge.
Figure~\ref{fig: 5} supports this for $\hat{\Gamma}_{1}$ of the $\operatorname{MA}(0)$
process. These non-symmetric matrices are under investigation.

%f3 #&#
%

%s4.2 #&#
\subsection{Applications} \label{sec: application}
%s4.2.1 #&#
\subsubsection{Order determination of a moving average process} \label
{subsubsec: maarorder} %Suppose we %are given a data set from a
%high-dimensional moving average process. Then d
A method to determine the order $q$ of a moving average process in the
univariate case
is to plot the correlogram (lag vs. sample autocorrelation graph) and
$\hat{q}$ is taken to be an estimate of $q$, if the sample
autocorrelations of order greater than $\hat{q}$ are small. %In
%high-dimensional case, as sample auotocovariances are not consistent
%for the corresponding population autocovariances, method of
%correleogram is useless.
In the high-dimensional case, as far as we know, there is no method in
the literature for estimating $q$.
%Here we use Remark~\ref{rem: diagnosis},

We use Theorem~\ref{thmm: lsd} to propose an analogous graphical method
of determining~$q$.
First, a look at Theorem~\ref{thmm: lsd} reveals that the LSD of $\hat
{\Gamma}_{u}\hat{\Gamma}_{u}^{*}$, for different $u$, can differ only
due to the distribution of
$\mathcal{C}_u=\{c_{j-j^{\prime}+u}: 0 \leq j,j^\prime\leq q\}$.
However, by applying (\ref{eqn: limitc}), it is not hard to see that
the joint distribution of
$\mathcal{C}_u$ are identical for all $u> q$ and are different for all
$0 \leq u \leq q$.

Therefore, when $X_t$ is a $\operatorname{MA}(q)$ process, the LSD of $\hat{\Gamma
}_{u}\hat{\Gamma}_{u}^{*}$ are identical for all $u > q$ and are
different for all $0 \leq u \leq q$.
% Consequently, under $\operatorname{MA}(\infty)$ process, the LSD of $\hat{\Gamma}_{u}
%\hat{\Gamma}_{u}^{*}$ are %different for all $u \geq0$.
These observations also hold true for any symmetric polynomial $\Pi_u$
in $\{\hat{\Gamma}_{u},\hat{\Gamma}_{u}^{*}\}$, for $u \geq0$.
%In Sections \ref{subsubsec: maarorder} and \ref{subsec: arorder}, we
%use these observations to suggest an appropriate order determination
%method. % for Moving Average process for a given data set. ???can we
%take a modfied version of thsi remark to the next section??
%\end{remark}

Let, for all $u \geq0$, $\Pi_u$ be a symmetric polynomial in $\{\hat
{\Gamma}_{u}, \hat\Gamma_u^*\}$. Note that the lower the order of the
polynomials, the lesser would be the moment conditions required for the
LSD to be valid. As an analogue of the correlogram, we propose to plot
the ECDF of some chosen $\Pi_u$ for first few sample autocovariance
matrices in the same graph. We say that $\hat{q}$ is an estimate of
$q$, if the ECDF of $\Pi_u$ with order $u > \hat{q}$ empirically
coincide with each other. %Note that, this method is valid only for
%large sample and simulations in previous section suggest that sample
%size $300$ is enough in this context. Also as weak convergence of ESD
%of $\Pi_i$ to its LSD is almost sure, it is easy to see that $\hat{q}$
%is a strong consistent estimator of $q$.
For example, consider the discussions in part (b) of Section~\ref
{subsec: example} and Figures \ref{fig: 1} and \ref{fig: 2}. There $q$
is determined quite accurately in the simulated data.
%As examples, we refer to Figures \ref{fig: 1} and \ref{fig: 2} and
%(b) ???? in Section~\ref{subsec: example}. %\textbf{????this IS
%Section~4!!, I dont understand}.
%\vskip10pt
%\begin{remark}
%One can also study the ESD of the asymmetric matrices $\{\hat{
%\Gamma}_{i}\}$, for $i \geq1$. Figure~\ref{fig: 5} supports the fact
%that ESD of $\hat{\Gamma}_{1}$ converges at least for $\operatorname{MA}(0)$ process.
%These need further investigation.
%s4.2.2 #&#
\subsubsection{Order determination of an autoregressive processes}
\label{subsec: arorder}
Another important problem is to determine the order of an infinite
dimensional vector Autoregressive (IVAR) process
%. Suppose we have data from an IVAR process
%
%e4.1 #&#
\begin{equation}
\label{eqn: ark} X_t = \varepsilon_t + A_1
X_{t-1} + A_2 X_{t-1} + \cdots+
A_{k}X_{t-k},
\end{equation}
where $k$ is unknown. Under suitable assumptions on the $p \times p$
parameter matrices $\{A_i\}$, one can show that $X_t$ satisfies (\ref
{eqn:model}) [see \citet{BB2013a}].
Suppose
% $\{X_t\}$ and $\{\varepsilon_t\}$ are all $p$-dimensional vectors.
$\{\varepsilon_t\}$ satisfies assumptions (A1) and (A2). %Suppose $
%\{A_i\}$ are $p \times p$ unknown parameter matrices.
Suppose the unknown parameter matrices $\{A_i\}$ are such that (\ref
{eqn: ark}) is stationary, and consistent estimators $\{\hat{A}_{i}\}$
for $\{A_i\}$ are available. By consistency, here we mean that the
limit of the spectral norm of $(\hat{A}_{i}-{A}_{i})$ is zero (in
probability). Such estimates are often available [see the end of
Section $3$ and discussions after Theorem $4.2$ in Section $4$ of \citet
{BB2013a}].
Also suppose that $\{A_i\}$ are compactly supported and for any finite
symmetric polynomial $\Pi$ of $\{A_i\}$
%\begin{equation}
$\lim p^{-1} E \operatorname{Tr} (\Pi) < \infty$ so that assumption
(A3) is
satisfied.

Then it is easy to see that, for each $u \geq0$, the LSD of $\hat
{\Gamma}_{u}\hat{\Gamma}_{u}^{*}$ for the process $\{\varepsilon_t\}$
[i.e., for the $\operatorname{MA}(0)$ process], coincides with the LSD (in
probability) for $\{\hat{\varepsilon}_{t}^{(k)}= X_t - \sum_{i=1}^{k}\hat{A}_{i}X_{t-i}\}$. See Section $8$ of the supplementary
file \citet{BB2014freesup} for the proof. Instead of $k$, if we use any
other positive integer $s < k$, then the residual process $\{\hat
{\varepsilon}_{t}^{(s)}\}$ does not behave like the $\operatorname{MA}(0)$ process.
% we show that, for a fix $u \geq0$, LSD of $\hat{\Gamma}_{u}\hat{
%\Gamma}_{u}^{*}$ for $\{\varepsilon_t\}$ and $\{\hat{
%\varepsilon}_{t}^{(k)}\}$ are identical, if $k$ is the order of the
%underlying IVAR process.
As ECDF of $\hat{\Gamma}_{u}\hat{\Gamma}_{u}^{*}$ for $u=1,2$ coincide
(almost surely) under $\operatorname{MA}(0)$ process, to determine the order of the
IVAR process, it is enough to check whether the ECDF of $\hat{\Gamma
}_{u}\hat{\Gamma}_{u}^{*}$ of $\{\hat{\varepsilon}_{t}^{(k)}\}$ for
$u=1,2$ coincides or not. Therefore, if we plot the ECDF of $\hat
{\Gamma
}_{u}\hat{\Gamma}_{u}^{*}$, $u=1,2$ for the residual process $\{\hat
{\varepsilon}_{t}^{(s)}\}$ in the same graph, the two distribution
functions coincide only when $s=k$. %then only the fitting of the
%autoregressive process is good and $k$ is the appropriate order of the
%underlying autoregressive process.
Hence, we may successively fit an $\operatorname{IVAR}(s)$ process for $s=0,1,2,\ldots
$ and for each $s$, plot the ECDF of $\hat{\Gamma}_{u}\hat{\Gamma
}_{u}^{*}$, $u=1,2$ for residuals $\{\hat{\varepsilon}_{t}^{(s)}\}$ in
the same graph. %and see whether two ECDFs are significantly close to
%each other or not. If they are close enough, then we say that $1$ is
%an estimate of order $k$ of the IVAR process. If the ECDFs are not
%close to each other, then continue the same process for $\{\hat{
%\varepsilon}_{t}^{(2)}\}$.
We say that $m$ is an estimate of order $k$ of the IVAR process, if the
ECDF of $\hat{\Gamma}_{u}\hat{\Gamma}_{u}^{*}$, $u=1,2$ do not coincide
for all $s <m$ and coincide for $s=m$.

For illustration, consider the following IVAR processes. Let
$\varepsilon_t \sim\mathcal{N}_{p}(0,I_p), \forall t$.

\textit{Model} $5$: $X_t = \varepsilon_t + 0.5 X_{t-1}$.

\textit{Model} $6$: $X_t = \varepsilon_t + 0.5 X_{t-1} + 0.2X_{t-2}$.

We let $p=n$ and draw a sample of size $n=500$. Assuming that
we do not know the parameter matrices, we use their banded estimator
from \citet{BB2013a}. % for the parameter matrices.

\begin{figure}

\includegraphics{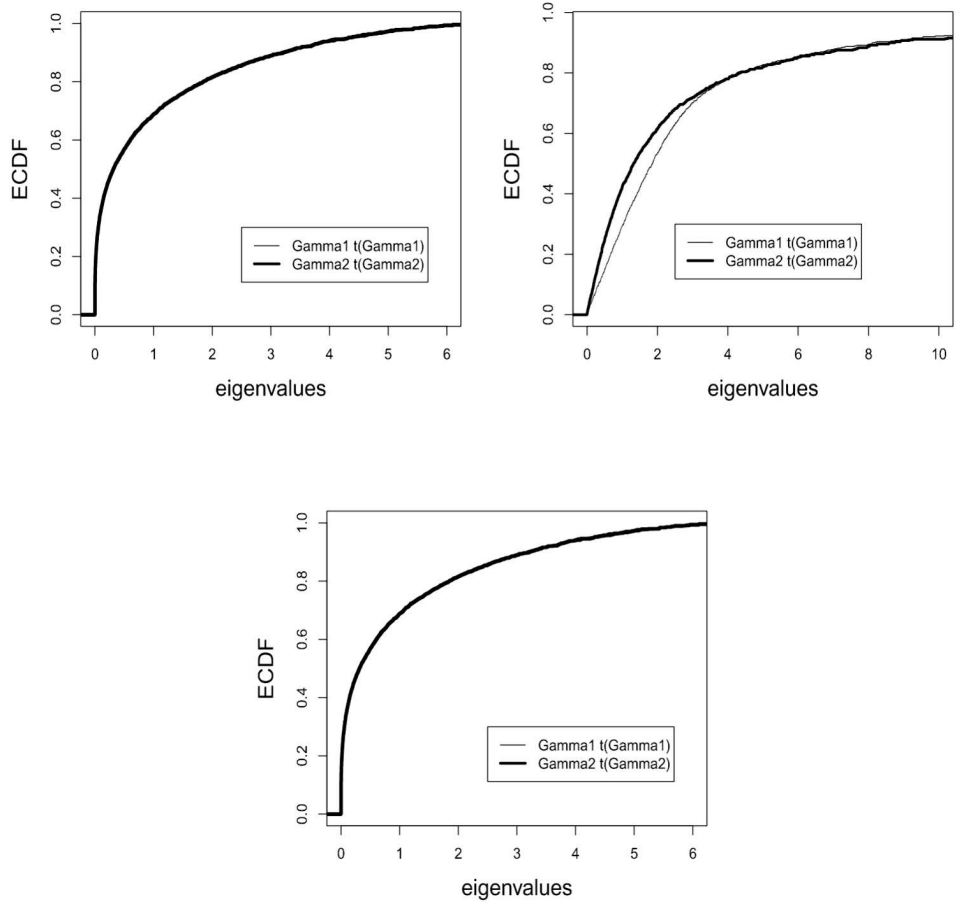}

\caption{$n=p=500$. Row $1$ left: ECDF of $\hat{\Gamma}_{1}\hat
{\Gamma}_{1}^{*}$ and $\hat
{\Gamma}_{2}\hat{\Gamma}_{2}^{*}$ for residuals after fitting $\operatorname{AR}(1)$
in model $5$. Row $1$ right: same after fitting $\operatorname{AR}(1)$ in model $6$.
Row $2$: same after fitting $\operatorname{AR}(2)$ in model~$6$.}
\label{fig: 3}
\end{figure}

For model $5$, we plot the two ECDFs of $\hat{\Gamma}_{u}\hat{\Gamma
}_{u}^{*}, u=1,2$ for the residual process $\{\hat{\varepsilon
}_{t}^{(1)}\}$ in the same graph and observe that they coincide. See
row $1$, left panel in Figure~\ref{fig: 3}. Therefore, $1$ is an
estimate of the order of model $5$. For model $6$, we do the same but
the two ECDFs do not coincide (see row $1$, right panel in Figure~\ref
{fig: 3}). In row $2$ of Figure~\ref{fig: 3}, the same two ECDFs are\vadjust{\goodbreak}
plotted for $\{\hat{\varepsilon}_{t}^{(2)}\}$ and they coincide and
hence $2$ is an estimate of the order for model $6$.

%f4 #&#

%s4.2.3 #&#
\subsubsection{Asymptotic distribution of traces and an application in
testing} \label{subsubsec: trace}
One of the referees raised the issue of convergence in distribution of
the trace of any autocovariance matrix and if such a result could be
possibly used for testing problems. %We show that this can be done by
%use of Lemma $4.1$ in the supplementary file. For all $T \geq1$, this
%lemma answers for $\lim E(\operatorname{Tr}(\pi)-E\operatorname{Tr}(
%\pi))^{T}$.
%As discussed at the beginning of this section, asymptotic distribution
%of appropriately scaled and centered spectral statistics can be found
%as a byproduct of Lemma $4.1$ in the supplementary file.
Let $\Pi:= \Pi(\hat{\Gamma}_{u},\hat{\Gamma}_{u}^{*}: u \geq0)$
be a
symmetric polynomial in $\{\hat{\Gamma}_{u},\hat{\Gamma}_{u}^{*}: u
\geq0\}$ and $\sigma_{\Pi}^{2} = \lim E(\operatorname{Tr}(\Pi
)-E\operatorname{Tr}(\Pi
))^{2} $. Then, for $d \geq1$, using some combinatorial calculations
[see Lemma $2.1$ in the supplementary file \citet{BB2014freesup}], we have
\begin{eqnarray*}
&&\lim E \bigl(\operatorname{Tr}(\Pi)-E\operatorname{Tr}(\Pi)
\bigr)^{T}
\\
&&\qquad = %
\cases{ 0,&\quad$\mbox{if $T=2d-1$,}$
\vspace*{2pt}\cr
\displaystyle\Biggl(\prod_{k=1}^{d}(2d-2k+1)
\Biggr) \sigma_{\Pi}^{2d},&\quad $\mbox{if $T=2d$}.$}
\end{eqnarray*}
Therefore,
%
%e4.2 #&#
\begin{equation}
\label{eqn: ad} \bigl(\operatorname{Tr}(\Pi)-E\operatorname{Tr}(\Pi) \bigr)
\stackrel {\mathcal{D}} {\rightarrow} \mathcal{N} \bigl(0,\sigma_{\Pi}^{2}
\bigr).
\end{equation}
The following are some examples and simulations to support
(\ref{eqn: ad}). We consider $n=p$ and $\varepsilon_t \sim\mathcal
{N}_{p}(0,I_{p})$, where $\varepsilon_t$'s are independent. % for
%Example~\ref{example: 5}. Similarly, all the detail calculations of
%Example~\ref{example: 4} can be done. Here we omit details for Example~\ref{example: 4}. %, but in can similarly be done using the same idea
%as in Example~\ref{example: 5}. %that the dimension $p$ equals the
%sample size $n$.

\begin{figure}

\includegraphics{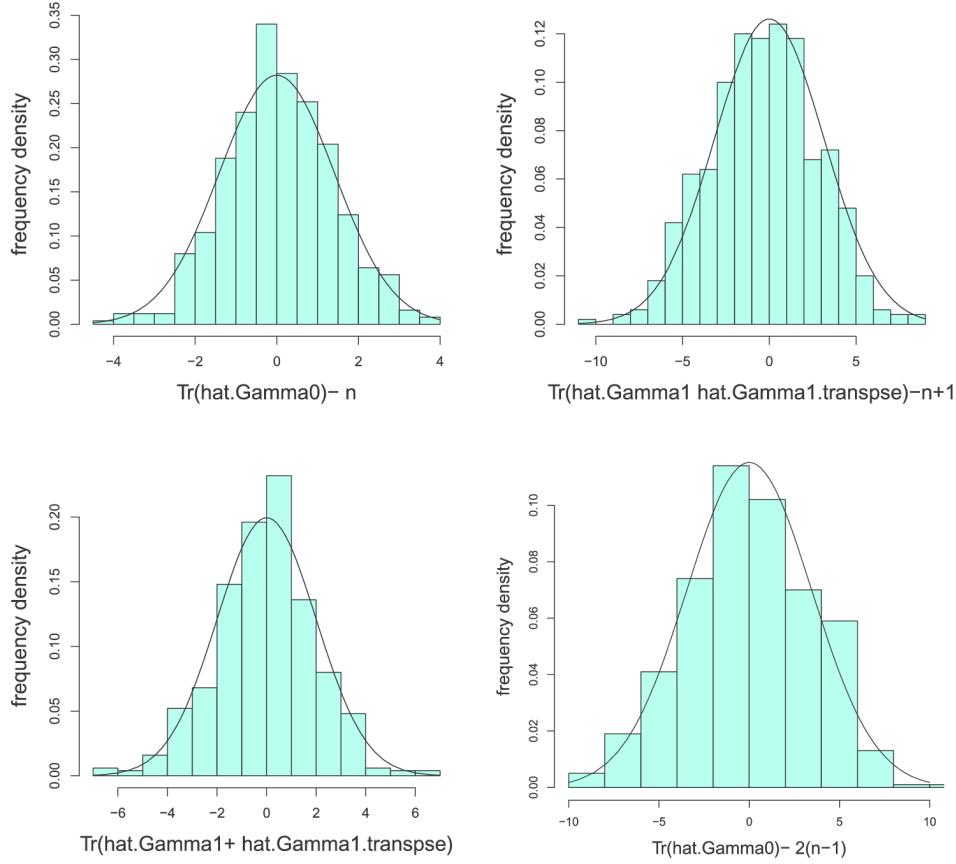}

\caption{$n=p=500$ and $500$ replications. Row (1)
left, rows (1) right and (2) left
represent respectively the histogram of $(\operatorname{Tr}(\hat
{\Gamma}_{0}) -
n)$, $(\operatorname{Tr}(\hat{\Gamma}_{1}\hat{\Gamma}_{1}^{*}) -
n+1)$ and
$\operatorname{Tr}(\hat{\Gamma}_{1}+\hat{\Gamma}_{1}^{*})$, when
$X_t =
\varepsilon_t$.
Row~(2) right represents the histogram of $(\operatorname{Tr}(\hat
{\Gamma
}_{0}) - 2n)$, when $X_t = \varepsilon_t + \varepsilon_{t-1}$.}\label
{fig: 4}
\end{figure}

%ex4 #&#
\begin{example} \label{example: 4}
Let $X_t = \varepsilon_t, \forall t$. Then $ E(\operatorname{Tr}\hat
{\Gamma
}_{0}) = n$, $E (\operatorname{Tr}(\hat{\Gamma}_{1}\hat{\Gamma
}_{1}^{*}) ) = n-1$, $E (\operatorname{Tr}(\hat{\Gamma
}_{1}+\hat{\Gamma
}_{1}^{*}) ) = 0$ and $\lim E(\operatorname{Tr}(\hat{\Gamma
}_{0})-E\operatorname
{Tr}(\hat{\Gamma}_{0}))^2=2$,
$\lim E (\operatorname{Tr}(\hat{\Gamma}_{1}\hat{\Gamma
}_{1}^{*}) -E\operatorname
{Tr}(\hat{\Gamma}_{1}\hat{\Gamma}_{1}^{*})  )^2= 10$, $\lim
E
(\operatorname{Tr}(\hat{\Gamma}_{1} + \hat{\Gamma}_{1}^{*}) -
E\operatorname{Tr}(\hat
{\Gamma}_{1}+\hat{\Gamma}_{1}^{*})  )^2 = 4$. We omit the detailed
calculations which are simpler than the calculations in the next
example. Hence,
%\begin{eqnarray} \label{eqn: clt1}
$(\operatorname{Tr}(\hat{\Gamma}_{0}) - n) \stackrel{\mathcal
{D}}{\rightarrow}
\mathcal{N}(0,2)$, %\nonumber\\
$(\operatorname{Tr}(\hat{\Gamma}_{1}\hat{\Gamma}_{1}^{*}) - n+1)
\stackrel
{\mathcal{D}}{\rightarrow} \mathcal{N}(0,10)$, and %\nonumber\\
$\operatorname{Tr}(\hat{\Gamma}_{1}+\hat{\Gamma}_{1}^{*})
\stackrel{\mathcal
{D}}{\rightarrow} \mathcal{N}(0,4)$. %\nonumber
%\end{eqnarray}
Simulation results given in rows $1$ and $2$, left panel,
Figure~\ref{fig: 4} support the above convergences.
%Note that $\operatorname{Tr}(\hat{\Gamma}_{0}) = n^{-1}\sum_{i,t}
%\varepsilon_{t,i}^{2}$ and hence one can also easily verify (\ref{eqn:
%clt1}) by the central limit theorem.
\end{example}

%ex5 #&#
\begin{example} \label{example: 5}
Let $X_t = \varepsilon_t + \varepsilon_{t-1}$. Then $ E(\operatorname
{Tr}(\hat
{\Gamma}_{0})) = 2(n-1)$, $\lim E(\operatorname{Tr}(\hat{\Gamma
}_{0})-E\operatorname
{Tr}(\hat{\Gamma}_{0}))^2=8$. In Section~\ref{subsec: detailtrace}, we
show details of this calculation. Hence,
%\begin{equation}
$ (\operatorname{Tr}(\hat{\Gamma}_{0}) - 2(n-1)) \stackrel{\mathcal
{D}}{\rightarrow} \mathcal{N}(0,8)$. %\nonumber
%\end{equation}
The simulation result given in row $2$, right panel, Figure~\ref{fig:
4} supports this convergence.
\end{example}

These results can be used for testing. For example, suppose
we wish to test
\[
H_0: X_t = \varepsilon_t\qquad \forall
t \quad\mbox{against}\quad H_1: X_t = \varepsilon_t +
\varepsilon_{t-1}\qquad \forall t.
\]
Then $(\operatorname{Tr}(\hat{\Gamma}_{0}) - n)$ can be used as a
test statistic and large value of the test statistic will imply
rejection of $H_0$.
Clearly, this idea can be extended to test other pairs of
%Similarly, one can perform many tests for high-dimensional linear
%process (\ref{eqn:model}) under
simple null and alternative hypotheses for model (\ref{eqn:model}).

%\section{Discussions} \label{sec: future}
%s5 #&#
\section{Proofs} \label{sec: proof}
We first prove Theorem~\ref{thmm: lsd}. For this purpose, we need the
following notions and results. % of convergence of sequence of NCPs
%which will provide us the theoretical basis for joint convergence of
%matrices.
%s5.1 #&#
\subsection{Convergence of NCPs and assymptotic freeness} \label
{subsec: ncpconvergence}
A sequence of NCPs $(\operatorname{Span}\{b_i^{(n)}\},\varphi_n)$ is
said to
converge to an NCP $(\operatorname{Span}\{b_i\},\varphi)$, if for any
polynomial~$\pi$,
%
%e5.1 #&#
\begin{equation}
\label{eqn: polyconvfree} \lim_{n} \varphi_{n} \bigl(\pi
\bigl(b_{i}^{(n)}: i \geq0 \bigr) \bigr) = \varphi \bigl(
\pi(b_{i}: i \geq0) \bigr).
\end{equation}
Suppose $(\operatorname{Span}\{b_{ij}^{(n)}: i \geq0, 1 \leq j \leq
k\},
\varphi_n)$ converges to $(\operatorname{Span}\{b_{ij}: i \geq0, 1
\leq j \leq
k\}, \varphi)$. Then $\operatorname{Span}\{b_{ij}^{(n)}: i \geq0\}$,
$1 \leq j
\leq k$, are said to be \textit{asymptotically free} if $\operatorname
{Span}\{
b_{ij}: i \geq0\}$ are free across $1 \leq j \leq k$.

%f5 #&#

%Let $\operatorname{Span}\{b_{ij}^{(n)}: i \geq0\}$, for all $1 \leq j
%\leq
%k$, be sequences of NCPs with same state $\varphi_n$. These $k$
%sequences of NCPs are said to be asymptotically free if, for each $j
%\geq0$, $(\operatorname{Span}\{b_{ij}^{(n)}: i \geq0\},\varphi_n)$
%converges
%to $(\operatorname{Span}\{b_{ij}: i \geq0\}, \varphi)$ (say) and these
%are
%free across $ 1\leq j \leq k$.

Let $W_{n\times n}$ be a \textit{Wigner} matrix. Let $\varphi_{n} =
n^{-1}E \operatorname{Tr}$. Let $\{B_{i,n}\}$ and $\{D_{i,n}\}$,
$1\leq i \leq
J$ be sequences of non-random, compactly supported, square matrices of
order $n$ each of which converges in the above sense. Then, under
assumption~(A2), the following facts are true. For (a) and (b),
see \citet{ZAG2010}. (c) follows from (a), (b) and Theorem $11.12$,
page $180$ of \citet{NS2006}. (d) is immediate from (a), (b) and
(c). We drop the suffix $n$ for clarity.
\begin{longlist}[(a)]
\item[(a)] $W/\sqrt{n}$ converges to the semi-circle law with
moment sequence (\ref{eqn: smoments}). % (see \citep{ZAG2010}).

\item[(b)] $W/\sqrt{n}$ and $\{B_{i}, D_{i}, 1\leq i\leq J\}$ are
asymptotically free. %(see \citep{ZAG2010}).

\item[(c)] $\{n^{-1}WB_iW, 1\leq i \leq J\}$ and $\{D_{i}, 1\leq
i\leq J\}$ are asymptotically free.

\item[(d)] To compute $\lim n^{-1} E\operatorname{Tr}(n^{-k}D_1 WB_1W D_2
WB_2W D_3\cdots WB_kW D_{k+1})$, one can assume that $W/\sqrt{n}$, $\{
B_i\}$ and $\{D_i\}$ are asymptotically free.
\end{longlist}

%s5.2 #&#
\subsection{Proof of Theorem \texorpdfstring{\protect\ref{thmm: lsd}}{3.1}} \label{subsec: lsdproof}
To prove Theorem~\ref{thmm: lsd}, as discussed in Section~\ref{sec:
Intro}, we have to show (M1), (M4) and (C) are satisfied. Here, we
shall only establish (M1). Proof of (M4) and (C) are given respectively
in Sections $3$ and $4$ of the supplementary file \citet{BB2014freesup}.

To establish (M1), we have to essentially show (\ref{eqn: lsdmoments}).
Let $\hat{\Gamma}_u (\varepsilon)$ be the $u$th order sample
autocovariance matrix of the process $\{ \varepsilon_t\}$. %: t=0,
%\pm1, \pm2, \cdots\}$.
Let
%
%e5.2 #&#
\begin{equation}
\label{eqn: Delta} \Delta_u = \sum_{j=0}^q
\sum_{j^\prime=0}^q \psi_j \hat{
\Gamma}_{j^\prime- j+u} (\varepsilon) \psi_{j^\prime}^{*}\qquad \forall
u \ge0.
\end{equation}
%
%Then for any polynomial $\Pi$, $\Pi(\hat{\Gamma}_{u},\hat{
%\Gamma}_{u}^{*}: u\geq0)$ and $\Pi(\Delta_{u},\Delta_{u}^{*}: u\geq
%0)$ have same limits in average trace. For proof see Section $1$ of
%Supplementary file \citet{BB2014freesup}. Therefore,
Then by Lemma $1.1$ of the supplementary file \citet{BB2014freesup}, it
is enough to show (\ref{eqn: lsdmoments}) after we replace $\{\hat
{\Gamma}_{u},\hat{\Gamma}_{u}^{*}\}$ by $\{\Delta_{u},\Delta
_{u}^{*}\}$.
%
%% with the same state $\varphi_{n}$.
%We define a $p \times n$ matrix
%%\begin{equation} \label{eqn: datamatrix}
%$Z=(\varepsilon_1  \varepsilon_2  \cdots\varepsilon_n)$.
%%\end{equation}
% It is
%called the ID matrix as all its elements are independently distributed
%with mean $0$
%and variance $1$. We embed $Z$ into a Wigner matrix $W$ of order
%$(n+p)$. Thus $$W= {
%W^{(1)}  \ Z \choose Z^{*} \ \ W^{(2)}},$$
%
% where $W^{(1)}$ and $W^{(2)}$ are two independent Wigner
%matrices of order $p$ and $n$ respectively and also independent of $Z$
%and satisfy Assumption (A2). For any $i$,
%let $P_i$
%be the $n \times n$ matrix with upper $i$th diagonal $1$. For any
%matrices $B$ and $D$ of order $p$ and $n$ respectively, let $\bar{B}$
%and $\underbar{D}$ of order $(n+p)$ be the matrices
%\begin{equation} \label{eqn: upperbar}\bar{B} = { B \ \ 0
%\choose0 \ \ 0},\ \ \underbar{D} = { 0 \ \ 0
%\choose0 \ \ D}. %\nonumber
%\end{equation} %For any matrix $D$ of order $n$, we
%%define $\underbar{D}$ of order $(n+p)$ as
%%\begin{equation} \label{eqn: underbar}
%%\underbar{D}= {0 \ \ 0 \choose0 \ \ D}.
%%\end{equation} %For all $i=0,1,2,3,\cdots, n-1$, we define $P_i$
%%to be the $n \times n$ matrix with upper $i$th diagonal $1$. Note
%%that $P_0=I_n$. For any matrix $D$ of order $n$, we
%%define $\underbar{D}$ of order $(n+p)$ as
%%\begin{equation} \label{eqn: underbar}
%%\underbar{D}= {0 \ \ 0 \choose0 \ \ D}.
%%\end{equation}
Now
%
%e5.3 #&#
\begin{equation}
\label{eqn: embedding} n (n+p)^{-1}\bar{\Delta}_i =
(n+p)^{-1}\sum_{j=0}^{q} \sum
_{j^\prime
=0}^{q} \bar{\psi}_j W
\underbar{P}_{j^\prime- j+i} W \bar{\psi}_{j^\prime}^{*}\qquad \forall i
\ge0.
\end{equation}
Note that by (a) of Section~\ref{subsec: ncpconvergence}, $W/\sqrt
{n+p}$ converges to the semi-circle law with moment sequence (\ref{eqn:
smoments}). Moreover, by (\ref{eqn: limitc}) and (\ref{eqn: limita}),
$\{\bar{\psi}_j\}$ and $\{\underbar{P}_{j}\}$ converge respectively to
$\{a_j\}$ and $\{c_j\}$. Also,
by (b), (c) and (d) of Section~\ref{subsec: ncpconvergence},
$s$, $\{a_j\}$ and $\{c_j\}$ are freely independent. Therefore, by
(\ref
{eqn: polyconvfree}),
(\ref{eqn: lsdmoments}) holds and (M1) is verified. Hence, proof of
Theorem~\ref{thmm: lsd} is complete. %with $\{\hat{\Gamma}_{u},\hat{
%\Gamma}_{u}^{*} \}$ replaced by $\{\hat{\Delta}_{u},\hat{
%\Delta}_{u}^{*}\}$ in the left side. Hence, by the approximation of $\{
%\hat{\Gamma}_{u},\hat{\Gamma}_{u}^{*} \}$ by $\{\hat{\Delta}_{u},\hat{
%\Delta}_{u}^{*}\}$ given in Section $3$ of the supplementary file,
%Condition $(M1)$ holds.

%Next we shall
%To justify
%Corollary~\ref{rem: final} and
%Example~\ref{rem: marchenko} (a), (b),. For this purpose,
Next, we need %to know computation of moments for free variables. In
%Section~\ref{subsec: algorithm}, we %discuss the
an algorithm for computing moments of a particular type of polynomials
of free variables. % which are useful This will be useful in justifying
%\ref{rem: marchecnko} (a), (b) and also derive a recursion formula
%for the moments of ???

%s5.3 #&#
\subsection{Algorithm to compute moments of free variables} \label
{subsec: algorithm}
%Although there is no closed form expression for arbitrary moments of
%the above LSD, we can calculate the moments recursively.
%using (\ref{eqn: freemoments}).
%Note that we need to get hold of
As we have discussed in Section~\ref{sec: free}, all joint moments of
free variables are computable in terms of the moments of the individual
variables. The algorithm for computing moments under freeness is
different from the product rule under usual independence.
Note that, for our purpose, a typical term in the moment calculations
[see, e.g., (\ref{eqn: 1ma})] is
%
%e5.4 #&#
\begin{equation}
\label{eqn: type} \varphi(d_{0}sb_{1}sd_{1}sb_{2}sd_{2}
\cdots sb_{n}sd_{n})\qquad\mbox{where $\{b_i\}$,
$ \{d_i\}$ and $s$ are free.}
\end{equation}
Note that in our case, since $\operatorname{Tr}(AB) = \operatorname
{Tr}(BA)$, our
$\varphi$ satisfies $\varphi(ab) = \varphi(ba),   \forall a,b$. In this
section, we shall discuss the algorithm for computing (\ref{eqn: type})
in terms of the moments of $\{b_i\}$, $\{d_i\}$ and $s$.

Let $\mathit{NC}(n)$ be the set of all non-crossing partitions of $\{1,2,\ldots,
n\}$. Define recursively a family of multiplicative, multi-linear
functionals $\varphi_{\pi}(n \geq1, \pi\in \mathit{NC}(n))$ by the following formula.
If $\pi= \{V_{1}, V_{2},\ldots, V_{r} \} \in \mathit{NC}(n)$, then
%
%e5.5 #&#
\begin{equation}\quad
\label{eqn: multi1} \varphi_{\pi}[a_{1},a_{2},
\ldots, a_{n}] := \varphi (V_1)[a_{1},a_{2},
\ldots,a_{n}]\cdots\varphi(V_r)[a_{1},a_{2},
\ldots,a_{n}],
\end{equation}
where
%
%e5.6 #&#
\begin{equation}\qquad
\label{eqn: multi2} \varphi(V)[a_{1},a_{2},\ldots,
a_{n}] := \varphi _{s}(a_{i_1}a_{i_2}
\cdots a_{i_s}) \qquad\mbox{for } V = (i_1 < i_2 <
\cdots< i_s).
\end{equation}
%
% $b_{1},b_{2},\ldots, b_{n} \in\mbox{Span}\{a_{i},a_{i}^{*}: i=0,1,2,
%\ldots\}$ and $d_{1},d_{2},\ldots, d_{n} \in\mbox{Span}
%\{c_{i},c_{i}^{*}: i=0,1,2,\ldots\}$.
Let $\mathit{NC}_{2}(2n)$ be the set of all non-crossing pair partitions of
$\{1,2,\ldots,2n\}$ and $K(\pi)\in \mathit{NC}(n)$ be the Kreweras
complement of
the partition $\pi$ [see Definition~$9.21$ in \citet
{NS2006}]. Then we have the following lemma. Relation (\ref{eqn:
freemoments}) will be useful to justify Example~\ref{rem: marchenko}(a)
and (b). %in Section~\ref{subsection: marchenko}
Relations (\ref{eqn: 11}) and (\ref{eqn: bb}) will be useful to prove
Lemma~\ref{rem: final}.

%le5.1 #&#
\begin{lemma} \label{eqn: momentalgo}\textup{(a)}
%
%e5.7 #&#
%e5.8 #&#
%e5.9 #&#
\begin{eqnarray}
\label{eqn: 11} && \varphi(d_{0}sb_{1}sd_{1}sb_{2}
\cdots sd_{n})
\nonumber
\\
&&\qquad= \sum_{\pi\in \mathit{NC}_{2}(2n)} \varphi_{K(\pi)}[
b_{1},d_{1},b_{2},d_{2},\ldots,
b_{n},d_{n}d_{0}]
\\
%&=& \sum_{\pi\in \mathit{NC}_{2}(2n)} \varphi_{K(\pi)}[b_{1},d_{1},b_{2},d_{2},
%\ldots,b_{n},d_{n}, d_{0}],
&&\qquad= \sum_{\pi\in \mathit{NC}(n)}
\varphi_{\pi}[b_{1},b_{2},\ldots,
b_{n}]\varphi _{K(\pi)}[d_{1},d_{2},
\ldots, d_{n}d_{0}] \label{eqn: 1i}
\\
&&\qquad = \sum_{\pi\in \mathit{NC}(n)} \varphi_{\pi
}[d_{1},d_{2},
\ldots, d_{n}d_{0}]\varphi_{K(\pi)}[b_{1},b_{2},
\ldots, b_{n}]. \label{eqn: freemoments}
\end{eqnarray}
\textup{(b)} Fix $1 = k_0 < k_1 < \cdots< k_t \leq n$ and the
following subset of $\mathit{NC}_{2}(2n)$ as
\[
\mathcal{S} = \bigl\{ \pi\in \mathit{NC}_{2}(2n): \{2k_{i},
2k_{i+1}-1\} \in\pi , 0 \leq i \leq t, k_{t+1} =
k_{0} \bigr\}.
\]
%
%\begin{lemma}
%(a) As $k_n(s)=1$ for $n=2$ and $0$ otherwise, by $(22.10)$ of
%\textnormal{\citet{NS2006}},
%have another expression for (\ref{eqn: type}) as
%(c) %Now let us get back (\ref{eqn: 11}).
%Consider any $\pi$ in $\mathit{NC}_{2}(2n)$ such that the blocks $\{1,2k_t\},
%\{2,2k_1-1\}, \{2k_1,$ $2k_2-1\},\ldots, \{2k_{t-1},2k_{t}-1\} \in
%\pi$ for any $1 < k_1 < k_2 < \ldots<k_t \leq n$. Then its Kreweras
%complement is $K(\pi)= \{\{\overline{2k_s-1}: 0 \leq s \leq t\}, \{
%\overline{2k_s},2k_s+1,\overline{2k_s+1},\ldots, 2k_{s+1}-2,
%\overline{2k_{s+1}-2}\}: 0 \leq s \leq t, k_0=1\}$. Note that the
%indices $\overline{2k_s-1}$ and $\overline{2k_s}$ are respectively
%occupied by $b_{k_s}$ and $d_{k_s}$. Hence, since $\varphi_{\pi}$ is
%multiplicative, we have %(\ref{eqn: type}) equals
Then
%
%e5.10 #&#
\begin{eqnarray}
\label{eqn: bb} && \sum_{\pi\in\mathcal{S}} \varphi_{K(\pi)}[b_1,
d_1, b_2, d_2,\ldots, b_n,
d_nd_{0}]
\nonumber
\\[-8pt]
\\[-8pt]
\nonumber
&&\qquad= \varphi \Biggl(\prod_{s=0}^{t}
b_{k_s} \Biggr) \prod_{s=1}^{t+1}
\varphi (d_{k_{s-1}}sb_{k_{s-1}+1}sd_{k_{s-1}+1}\cdots s
d_{k_s-1}),
\end{eqnarray}
where $k_{0} = 1, d_{k_{t+1}-1} = d_nd_0$.
\end{lemma}

Relation (\ref{eqn: 11}) follows by (22.10) of
\citet{NS2006}. By freeness of $\{b_i\}$ and $\{d_i\}$, and by
properties of the Kreweras complement [Exercises 9.41(1), 9.42(1) and
(2) in \citet{NS2006}], (\ref{eqn: 1i}) and
(\ref{eqn: freemoments}) follow from (\ref{eqn: 11}). Relation (\ref
{eqn: bb}) follows from the multiplicative property (\ref{eqn: multi1})
and (\ref{eqn: multi2}) of partitions and from certain properties of
Kreweras complement. A detailed proof of (\ref{eqn: bb}) is given in
Section $9$ of the supplementary file \citet{BB2014freesup}.
\subsection{A recursion formula for moments and its proof} \label
{subsection: final}
In this section, we shall prove a lemma that provides a recursion
formula for the moments of the LSD of $2^{-1}(\hat{\Gamma}_{u} + \hat
{\Gamma}_{u}^{*})$, which will be used in the proof of Remark~\ref{rem:
liufinal} in the next section.
%This formula is cumbersome but is quite useful to us.
% (see ???)

Let
\[
D= 2^{-1} \bigl(\hat{\Gamma}_{u}+\hat{\Gamma}_{u}^{*}
\bigr)\quad\mbox{and}\quad d_{uq}=2^{-1} \bigl(\gamma_{uq}+
\gamma_{uq}^{*} \bigr),
\]
where $\gamma_{uq}$ is as in (\ref{eqn: gammafree}). Suppose $\theta$
is a $U(0,2\pi)$ variable, which is (classical) independent and
commutative with $\{{a}_j\}$ and ${d}_{uq}$. Recall $h(\lambda, \theta
)$ in (\ref{eqn: hhr}). %$\nu$ is commutative with $\{a_l\}$ and $\bar{
%\gamma}_{uq}$.
For any polynomial $\Pi= \Pi(\psi_j,\psi_j^{*}: j \geq0)$, let
$\bar
{\Pi} =\Pi(a_j,a_j^{*}: j \geq0)$.
%\begin{eqnarray}
%\label{eqn: h} h(\lambda,\theta) = (\sum_{j=0}^{\infty}e^{ij
%\theta}a_{j})(\sum_{j=0}^{\infty}e^{-ij\theta}a_{j}^{*}),\ \lambda=
%\{a_j, a_j^{*}: j \geq0\}
%\end{eqnarray}
For all $j \geq0$, let
%
%e5.11 #&#
%e5.12 #&#
\begin{eqnarray}
\label{eqn: rdef} \label{eqn: r} R_{uj}(\theta) &=& y^{-1}
(1+y)\varphi \bigl(h(\lambda ,\theta){d}_{uq}^{j-1} \bigr),
\\
\label{eqn: s} S_{uj}(\theta, \Pi) &=& y^{-1} (1+y)\varphi
\bigl(\bar{\Pi }h(\lambda,\theta){d}_{uq}^{j-1} \bigr).
\end{eqnarray}
\begin{lemma} \label{rem: final}
Let $X_t \sim\operatorname{MA}(q)$ process and suppose assumptions \textup{(A1)},
\textup{(A2)}, \textup{(A3)} and $pn^{-1} \to y \in(0,\infty)$ hold. Then for any
polynomial $\Pi= \Pi(\psi_j, \psi_j^{*}: j \geq0)$, we have
%
%e5.13 #&#
\begin{eqnarray}
\label{eqn: moment01} && \lim p^{-1} E \operatorname{Tr} \bigl(\Pi
D^{r} \bigr) %\nonumber
\nonumber
\\[-8pt]
\\[-8pt]
\nonumber
&&\qquad= \frac{1}{y}\sum_{t=1}^{r}
E_{\theta} \Biggl[ \bigl(y\cos(u\theta ) \bigr)^{t}\mathop{ \sum
_{1 \leq i_1, i_2,\ldots, i_t \leq r}}_{\sum
_{j=1}^{t}i_{j}=r} S_{ui_1}(\theta, \Pi) \Biggl(\prod
_{k=2}^{t}R_{u
i_k}(\theta) \Biggr)
\Biggr].
%\bigg(\sum_{\pi\in\mathcal{P}_{r}} \bigg(\frac{t!}{
%\prod_{k=1}^{s}f_{k}!}y^{t-1}(\cos(i\nu))^{t}\prod_{k=1}^{t}R_{u i_k}(
%\nu)\bigg)\bigg)\bigg].
\end{eqnarray}
\end{lemma}
%
%\textbf{???comment about proof???}

\begin{pf}
From the proof of (\ref{eqn: lsdmoments}),
%Theorem~\ref{thmm: lsd},
it is immediate that
%
%e5.14 #&#
\begin{eqnarray}
\label{eqn: 1ststep}&&  \frac{y}{1+y}\lim p^{-1} E \operatorname{Tr}
\bigl(\Pi D^{r} \bigr)\nonumber \\
&&\qquad= \varphi \bigl(\bar{\Pi} {d}_{uq}^{r}
\bigr)
\nonumber
\\[-8pt]
\\[-8pt]
\nonumber
%\nonumber\\
%\nonumber\\ %\ \ \operatorname{by (\ref{eqn: jointlsd})} \nonumber\\
&&\qquad= \mathop{\sum_{j_k,j_k^\prime=1}}_{1\leq k \leq r}^{q}
\mathop{\sum_{v_k = u,-u}}_{1 \leq k \leq r} \varphi \Biggl(\bar{\Pi }\prod
_{k=1}^{r} {a}_{j_k}sc_{j_k-j_k^\prime+ v_k}s
{a}_{j_k^\prime
}^{*} \Biggr)\qquad\mbox{by (\ref{eqn: gammafree})}
\\
&&\qquad= \sum_{\sigma\in \mathit{NC}_{2}(2r)} \tau_{\sigma} \qquad\mbox{by
(\ref{eqn:
11})},\nonumber %\nonumber%\varphi_{\sigma}[c_{j_1-j_1^\prime+v_1}, a_{j_1^
%\prime}a_{j_2}, c_{j_2-j_2^\prime+ v_2} ]
\end{eqnarray}
where
\[
\tau_{\sigma} = \mathop{\sum_{j_k,j_k^\prime=1}}_{1\leq k \leq r}^{q}
\mathop{\sum_{v_k = u,-u}}_{{1 \leq k \leq r}}\varphi_{K(\sigma
)}
\bigl[c_{j_1-j_1^\prime+v_1}, {a}_{j_1^\prime}^{*} {a}_{j_2},
c_{j_2-j_2^\prime+ v_2}, {a}_{j_2^\prime}^{*} {a}_{j_3},\ldots,
{a}_{j_r^\prime}^{*} {a}_{j_1}\bar{\Pi} \bigr],
\]
and $K(\sigma)$ is the Kreweras complement of $\sigma$.
Now to compute (\ref{eqn: 1ststep}), we consider the decomposition of
$\mathit{NC}_{2}(2r) = \bigcup_{t=1}^{r}\mathcal{P}_{t}^{2r}$, where $\mathcal
{P}_{1}^{2r} = \{\sigma\in \mathit{NC}_{2}(2r): \{1,2\} \in\sigma\}$ and for
all $2 \leq t \leq r$,
\begin{eqnarray*}
\mathcal{P}_{t}^{2r} &=& \bigl\{\sigma\in
\mathit{NC}_{2}(2r): \{2k_0-1,2k_t\},\{
2k_0,2k_1-1\},\{2k_1,2k_2-1\},
\ldots,
\\
&& \{2k_{t-2},2k_{t-1}-1\} \in\sigma, 1 = k_0
<k_1 <k_2 <\cdots<k_{t-1} \leq r \bigr\}.
\end{eqnarray*}
Hence, (\ref{eqn: 1ststep}) is equivalent to
%
%e5.15 #&#
\begin{equation}
\label{eqn: 111} y(1+y)^{-1} \lim p^{-1}E \operatorname{Tr}
\bigl(\Pi D^{r} \bigr)= \sum_{t=1}^{r}
\mathcal{T}_{t},
\end{equation}
where for all $1 \leq t \leq r$,
%
%e5.16 #&#
\begin{equation}
\label{eqn: 2ndstep} \mathcal{T}_t = \sum_{\sigma\in\mathcal{P}_{t}^{2r}}
\tau _{\sigma} = \sum_{1=k_0<k_1<k_2< \cdots<k_{t-1} \leq r} g(t,
k_1,k_2,\ldots, k_{t-1}), %\nonumber
\end{equation}
and %$2^{t+1}(1+y)^{-t-1}g(t+1,k_1,k_2,\ldots, k_t)$ equals (we take
%$j_{k_{t+1}} = j_{k_{0}}$, $\bar{a}_{j_{k_{t+1}}}=\bar{a}_{j_{k_0}}
%\bar{\Pi}$, ${a}_{j_{k_{t+1}}}= {a}_{j_{k_0}}{\Pi}$, $C_{y} =
%\frac{y^{t+1}}{(1+y)^{t+2}}$)
%
%e5.17 #&#
%e5.18 #&#
\begin{eqnarray}
\label{eqn: 33} && 2^{t+1}(1+y)^{-t-1}g(t+1,k_1,k_2,
\ldots, k_t)\nonumber
\\
&&\qquad=\mathop{\sum_{ 1 \leq j_{k_s},j_{k_{s}}^\prime\leq q}}_
{v_{k_s}=u,-u} \varphi \Biggl(\prod
_{s=0}^{t} c_{{j_{k_s}}-j_{k_{s}}^\prime+
v_{k_s}} \Biggr) \prod
_{s=0}^{t} \varphi \bigl({a}_{j_{k_s}^\prime
}^{*}{d}_{uq}^{k_{s+1}-k_{s}-1}
{a}_{j_{k_{(s+1)}}} \bigr)
\\
&&\qquad \eqntext{\bigl[\mbox{by (\ref{eqn: bb}) and where $k_{t+1}=r+1$ and
${a}_{j_{k_{t+1}}}= {a}_{j_{k_0}}\bar{\Pi}$}\bigr].}
\end{eqnarray}
Now, by (\ref{eqn: limitc}),
%
%e5.19 #&#
\begin{eqnarray}
\varphi \Biggl(\prod_{s=0}^{t}
c_{{j_{k_s}}-j_{k_{s}}^\prime+ v_{k_s}} \Biggr) &=& %
\cases{ 1,&\quad $\mbox{if $\displaystyle\sum
_{s=0}^{t} \bigl({j_{k_s}}-j_{k_{s}}^\prime+
v_{k_s} \bigr)=0$},$ \vspace*{2pt}
\cr
0,& \quad$\mbox{otherwise}$}
\nonumber\\
&=& E_{\theta} \bigl(e^{i \theta\sum_{s=0}^{t} ({j_{k_s}}-j_{k_{s}}^\prime
+ v_{k_s})} \bigr)\nonumber\\
\eqntext{\mbox{where $\theta\sim
U(0, 2\pi)$}.}
\end{eqnarray}
Therefore, (\ref{eqn: 33}) is equivalent to
%
%e5.20 #&#
\begin{eqnarray}
&=& \frac{1}{1+y} \mathop{\sum_{ 1 \leq
j_{k_s},j_{k_{s}}^\prime\leq q} }_{{v_{k_s}=u,-u}} E_{\theta}
\bigl(e^{\sum
_{s=0}^{t}i \theta(j_{k_s} - j_{k_s}^\prime+ v_{k_s})} \bigr) \prod_{s=0}^{t}
\varphi \bigl({a}_{j_{k_s}^\prime
}^{*}{d}_{uq}^{k_{s+1}-k_{s}-1}
{a}_{j_{k_{(s+1)}}} \bigr)\nonumber
\\
%&=& C_{y} E_{\theta}\bigg(\sum_{ \stackrel{1 \leq j_{k_s},j_{k_{s}}^
%\prime\leq q} {v_{k_s}=u,-u}}
%\prod_{s=0}^{t}(e^{i \theta(j_{k_s} - j_{k_s}^\prime+ v_{k_s})})
%\varphi\big(a_{j_{k_s}^\prime}^{*} {d}_{uq}^{k_{s+1}-k_{s}-1}
%a_{j_{k_{(s+1)}}}\big)\bigg) \nonumber\\
&=&
\frac{1}{1+y} E_{\theta} \Biggl(\prod_{s=0}^{t}
\sum_{j_{k_s}^\prime,
j_{k_{s+1}},v_{k_s}} e^{i\theta(-j_{k_s}^\prime+
j_{k_{(s+1)}}+v_{k_s})}\varphi
\bigl(a_{j_{k_s}^\prime}^{*} {d}_{uq}^{k_{s+1}-k_{s}-1}
a_{j_{k_{(s+1)}}} \bigr) \Biggr)\nonumber
\\
\eqntext{\mbox{where $j_{k_{t+1}} = j_{k_0}$}.}
\end{eqnarray}
Note that for a fix $0 \leq s \leq t$, since $v_{k_s} = u, -u$,
\begin{eqnarray*}
\label{eqn: 3rdstep} && \sum_{j_{k_s}^\prime, j_{k_{s+1}},v_{k_s}} e^{i\theta
(-j_{k_s}^\prime+ j_{k_{(s+1)}}+v_{k_s})}
\varphi \bigl(a_{j_{k_s}^\prime
}^{*}{d}_{uq}^{k_{s+1}-k_{s}-1}
a_{j_{k_{(s+1)}}} \bigr)
\\
&&\qquad= 2\cos(u\theta) \varphi \biggl( \biggl(\sum_{}
e^{i\theta
j_{k_{(s+1)}}}a_{j_{k_{(s+1)}}} \biggr) \biggl(\sum
_{j_{k_s}}a_{j_{k_s}^\prime}^{*} e^{-i \theta j_{k_s}^\prime}
\biggr){d}_{uq}^{k_{s+1}-k_{s}-1} \biggr)
\nonumber
\\
%&=& 2\cos(u\theta) \varphi\bigg(h(\lambda,\theta)
%{d}_{uq}^{k_{s+1}-k_{s}-1}\bigg) \nonumber\\
&&\qquad= %
\cases{\displaystyle 2\cos(u\theta)
\frac{y}{1+y}R_{u(k_{s+1}-k_{s})}(\theta ),&\quad $1 \leq s \leq t-1,$ \vspace*{2pt}
\cr
2
\cos(u\theta)\displaystyle\frac{y}{1+y} S_{u(r+1-k_{t})}(\theta, \Pi),&\quad $s=t.$}
\end{eqnarray*}
Hence, for all $1=k_0<k_1<k_2< \cdots< k_{t} \leq r$, we have
\begin{eqnarray*}
&& g(t+1,k_1,k_2,\ldots, k_t)
\\
&&\qquad= y^{t+1}(1+y)^{-1} E_{\theta} \Biggl(\cos(u\theta
)S_{u(r+1-k_{t})}(\theta, \Pi)\prod_{s=0}^{t-1}
\cos(u \theta) R_{u(k_{s+1}-k_{s})}(\theta) \Biggr).
\end{eqnarray*}
Therefore, by (\ref{eqn: 2ndstep}), for all $1 \leq t \leq r$,
\begin{eqnarray*}
\mathcal{T}_{t} &=& \frac{y^{t}}{(1+y)} \mathop{\sum
_{1=k_0<k_1<
\cdots< k_{t} \leq r}}_{{k_{t+1}=r+1}}E_{\theta} \Biggl[ \bigl(\cos(u\theta )
\bigr)^{t}S_{u(r+1-k_{t})}(\theta, \Pi)\prod
_{s=0}^{t-1}R_{u(k_{s+1}-k_{s})}(\theta) \Biggr]
\\
&=& y^{t}(1+y)^{-1}\mathop{\sum_{1 \leq i_1, i_2 ,\ldots, i_t \leq
r}}_{i_1+i_2+\cdots+ i_t = r}
E_{\theta} \Biggl( \bigl(\cos(u\theta) \bigr)^{t}S_{u
i_{1}}(
\theta, \Pi)\prod_{s=2}^{t}
R_{ui_s}(\theta) \Biggr).
\end{eqnarray*}
Hence, by using (\ref{eqn: 111}) and (\ref{eqn: 2ndstep}), relation
(\ref{eqn: moment01}) follows.
%\begin{equation}
%\lim p^{-1}E\operatorname{Tr}(\Pi D^{r}) = y^{-1}\sum_{t=1}^{r} E_{
%\theta}
%\bigg[ (y\cos(u \theta))^t \bigg(\sum_{\stackrel{1 \leq i_1, i_2 ,
%\ldots, i_t \leq r}{i_1+i_2+\ldots+ i_t = r}} S_{u i_{1}}(\theta, \Pi)
%\prod_{s=2}^{t}R_{ui_s}(\theta))\bigg)\bigg] \nonumber
%\end{equation}
%and the proof of Lemma~\ref{rem: final} is complete.
\end{pf}

%s5.5 #&#
\subsection{Proof of Remark \texorpdfstring{\protect\ref{rem: liufinal}}{3.2}} \label
{subsection: liufinal}
(a) We prove (\ref{eqn: mmr})--(\ref{eqn: bbr}) stated in Remark~\ref
{rem: liufinal} under assumptions (A1)--(A3). By Remark~\ref
{remark: truncation}, assumption (A2) can be replaced by assumptions
(A4) and (A5).

Let us define
\begin{eqnarray*}
K(z,\theta) &=& \sum_{i=1}^{\infty}z^{-i}R_{ui}(
\theta),\qquad D = \sum_{i=1}^{\infty}z^{-i}
d_{uq}^{i},
\\
B(\lambda,z) &=& E_{\theta} \bigl( \cos(u\theta) h(\lambda,\theta )
\bigl(1+ y \cos(u \theta) K(z,\theta) \bigr)^{-1} \bigr),
\end{eqnarray*}
where $\{R_{ui}\}$ is defined in (\ref{eqn: rdef}). Now note that
%
%e5.21 #&#
\begin{eqnarray}
\label{eqn: yk} y(1+y)^{-1}K(z,\theta) &=& \sum
_{i=1}^{\infty} z^{-i} \varphi \bigl(h(\lambda
,\theta)d_{uq}^{i-1} \bigr)\qquad \mbox{by (\ref{eqn: r})}
\nonumber
\\[-8pt]
\\[-8pt]
\nonumber
&=& - z^{-1} \varphi \bigl(h(\lambda, \theta) \bigr) -
z^{-1}\varphi \bigl(h(\lambda ,\theta)D \bigr),
\end{eqnarray}
where
%
%e5.22 #&#
\begin{eqnarray}
&& \varphi \bigl(h(\lambda,\theta)D \bigr)\nonumber
\\
&&\qquad= \sum_{r=1}^{\infty} z^{-r}
y^{-1}\sum_{t=1}^{r}
E_{\theta
^\prime
} \Biggl[ \bigl(y\cos \bigl(u \theta^\prime \bigr)
\bigr)^t \nonumber\\
&&\qquad\quad{}\times \Biggl(\mathop{\sum_{1 \leq i_1,
i_2 ,\ldots, i_t \leq r}}_{{i_1+i_2+\cdots+ i_t = r}}
S_{u i_{1}} \bigl(\theta ^\prime, h(\lambda,\theta) \bigr)\prod
_{s=2}^{t}R_{ui_s} \bigl(
\theta^\prime \bigr) \Biggr)  \Biggr]\qquad \mbox{(by Lemma~\ref{rem: final})}\nonumber
\\
&&\qquad= y^{-1}\sum_{t=1}^{\infty}
E_{\theta^\prime} \Biggl[ \bigl(y\cos \bigl(u \theta ^\prime \bigr)
\bigr)^t \nonumber\\
&&\qquad\quad{}\times \sum_{r=t}^{\infty}
z^{-r} \Biggl(\mathop{\sum_{1
\leq
i_1, i_2 ,\ldots, i_t \leq r}}_{i_1+i_2+\cdots+ i_t = r} S_{u
i_{1}}
\bigl(\theta^\prime, h(\lambda,\theta) \bigr)\prod
_{s=2}^{t}R_{ui_s} \bigl(\theta
^\prime \bigr) \Biggr)  \Biggr]
\nonumber
\\
&&\qquad= y^{-1}\sum_{t=1}^{\infty}
E_{\theta^\prime} \Biggl[ \bigl(y\cos \bigl(u \theta ^\prime \bigr)
\bigr)^t \bigl(K \bigl(z,\theta^\prime \bigr)
\bigr)^{t-1} \Biggl( \sum_{r=1}^{\infty}
z^{-r} S_{u r} \bigl(\theta^\prime, h(\lambda,\theta)
\bigr) \Biggr) \Biggr]\nonumber
\\
&&\qquad= E_{\theta^\prime} \Biggl[ \cos \bigl(u\theta^\prime \bigr) \bigl(1+ y
\cos \bigl(u \theta ^\prime \bigr) K \bigl(z,\theta^\prime \bigr)
\bigr)^{-1} \Biggl( \sum_{r=1}^{\infty}
z^{-r} S_{u
r} \bigl(\theta^\prime, h(\lambda,\theta)
\bigr) \Biggr) \Biggr]\nonumber
\\
&&\qquad= \varphi \Biggl(\sum_{r=1}^{\infty}z^{-r}
d_{uq}^{r-1} h(\lambda ,\theta ) E_{\theta^\prime} \bigl[ h
\bigl(\lambda,\theta^\prime \bigr)\cos \bigl(u\theta ^\prime
\bigr) \bigl(1+ y \cos \bigl(u \theta^\prime \bigr) K \bigl(z,
\theta^\prime \bigr) \bigr)^{-1} \bigr] \Biggr)\nonumber
\\
&&\qquad= z^{-1}\varphi \bigl(h(\lambda,\theta)B(\lambda,z) \bigr) +
z^{-1} \varphi \bigl(Dh(\lambda,\theta)B(\lambda,z) \bigr).\nonumber
%\stackrel{\operatorname{similarly}}{=}& z^{-1}\varphi(h(\lambda,
%\theta)B(
%\lambda,z)) + z^{-2}\varphi(h(\lambda,\theta)B^{2}(\lambda,z))
% z^{-2} \varphi(Dh(\lambda,\theta)B^{2}(\lambda,z)).\nonumber
\end{eqnarray}
In a similar fashion, using $h(\lambda, \theta)B(\lambda,z)$
instead of $h(\lambda, \theta)$ in the above steps,
\[
\varphi \bigl(Dh(\lambda,\theta)B(\lambda,z) \bigr) = z^{-1}\varphi
\bigl(h(\lambda ,\theta )B^{2}(\lambda,z) \bigr) + z^{-1}
\varphi \bigl(Dh(\lambda,\theta)B^{2}(\lambda ,z) \bigr).
\]
Finally iterating we have
%
%e5.23 #&#
\begin{equation}
\label{eqn: finaliteration} \varphi \bigl(h(\lambda,\theta)D \bigr) = \sum
_{r=1}^{\infty} z^{-r}\varphi \bigl(h(\lambda,
\theta)B^{r}(\lambda,z) \bigr).
\end{equation}
We now need to show only (\ref{eqn: mmr}) and (\ref{eqn:
kkr}). Using (\ref{eqn: finaliteration}) and (\ref{eqn: yk}),
%noindent Proceeding in this way, we have
%
%e5.24 #&#
\begin{eqnarray}
\label{eqn: k} &&y(1+y)^{-1} K(z,\theta)
\nonumber
\\[-8pt]
\\[-8pt]
\nonumber
&&\qquad= - z^{-1} \varphi \Biggl(\sum_{r=0}^{\infty}
h(\lambda, \theta) z^{-r} B^{r}(\lambda,z) \Biggr) = \varphi
\bigl(h(\lambda, \theta) \bigl(B(\lambda,z)-z \bigr)^{-1} \bigr),
\end{eqnarray}
which is (\ref{eqn: kkr}) in Remark~\ref{rem: liufinal}. %Note that
%under (B), (\ref{eqn: kkr}) reduces to (\ref{eqn: kr}). %Let $F_{a}$
%be the joint distribution function of $\{\psi_{j},\psi_{j}^{*}: j\geq0
%\}$ ????. Then
%\begin{eqnarray} \label{eqn:k}
% K(z,\theta) %\nonumber\\
%&=& \int\bigg( E_{\theta^\prime}\big(\frac{\cos(u\theta^\prime) (h(
%\lambda,\theta^\prime) }{1+ y \cos(u \theta^\prime) K(z,\theta^\prime)}
%\big)-z\bigg)^{-1} h(\lambda, \theta) dF_{a}. %\nonumber
%\end{eqnarray}

Note that the above steps from (\ref{eqn: yk}) leading to (\ref{eqn:
k}) remain valid if we replace $h(\lambda,\theta)$ by $1$ in (\ref{eqn:
yk}). This yields [instead of (\ref{eqn: k})],
%
%e5.25 #&#
\begin{equation}
\label{eqn: m} y(1+y)^{-1}m(z) %\nonumber\\
= \varphi \bigl( \bigl(B(
\lambda,z)-z \bigr)^{-1} \bigr), %\int\bigg( E_{\theta^\prime}\big(
%\frac{\cos(u\theta^\prime) (h(\lambda,\theta^\prime) }{1+ y \cos(u
%\theta^\prime) K(z,\theta^\prime)}\big)-z\bigg)^{-1} dF_{a}, %\nonumber
\end{equation}
which is (\ref{eqn: mmr}) in Remark~\ref{rem: liufinal}.
Hence, the proof of Remark~\ref{rem: liufinal}(a) is complete. %It is easy
%to see that under (B), (\ref{eqn: mmr})-(\ref{eqn: bbr}) reduce to (
%\ref{eqn: mr})-(\ref{eqn: hr}).

(b) We now prove that under assumption (B), the Stieltjes
transform equations (\ref{eqn: mmr})--(\ref{eqn: bbr}) reduce to
equations (\ref{eqn: mr})--(\ref{eqn: hr}).

Note that under assumption (B),
%
%e5.26 #&#
\begin{eqnarray}
\label{eqn: newh} h(\lambda,\theta) &=& \Biggl(\sum_{j=0}^{\infty}e^{ij\theta}a_{j}
\Biggr) \Biggl(\sum_{j=0}^{\infty}e^{-ij\theta}a_{j}^{*}
\Biggr), \qquad\lambda= \bigl\{a_j, a_j^{*}: j \geq0
\bigr\}
\nonumber
\\
&=& y(1+y)^{-1} \Biggl(\sum_{j=0}^{\infty}e^{ij\theta}f_{j}(
\alpha) \Biggr) \Biggl(\sum_{j=0}^{\infty}e^{-ij\theta}f_{j}^{*}(
\alpha) \Biggr) \nonumber\\
&&{}+ (1+y)^{-1}\delta _{0}, \qquad\alpha\sim
F_{a}
\\
&=& y(1+y)^{-1}\Biggl|\sum_{j=0}^{\infty}e^{ij\theta}f_{j}(
\alpha)\Biggr|^{2} + (1+y)^{-1}\delta_{0}
\nonumber
\\
&=& y(1+y)^{-1}h_{1}(\alpha, \theta)+ (1+y)^{-1}
\delta_{0}, \qquad\alpha \sim F_{a},\nonumber
\end{eqnarray}
where $\delta_0$ is the degenerate random variable at $0$. Therefore,
after substituting (\ref{eqn: bbr}) in (\ref{eqn: kkr}), we have
\begin{eqnarray*}
K(z,\theta) & = &y^{-1}(1+y) \varphi \bigl(h(\lambda, \theta) \bigl(B(
\lambda ,z)-z \bigr)^{-1} \bigr)
\\
&=& y^{-1}(1+y) \varphi \biggl( \biggl[E_{\theta^\prime} \biggl(
\frac
{\cos
(u\theta^\prime) h(\lambda,\theta^\prime)}{ (1+ y \cos(u \theta
^\prime)
K(z,\theta^\prime))} \biggr) -z \biggr]^{-1}h(\lambda,\theta) \biggr)
\\
&=& \int \biggl[E_{\theta^\prime} \biggl( \frac{\cos(u\theta^\prime)
h_{1}(\alpha,\theta^\prime)}{ (1+ y \cos(u \theta^\prime)
K(z,\theta
^\prime))} \biggr) -z
\biggr]^{-1}h_{1}(\alpha,\theta) \,dF_{a}\qquad
\mbox{by (\ref{eqn: newh})}.
\end{eqnarray*}
Hence, (\ref{eqn: kkr}) reduces to (\ref{eqn: kr}).
Similarly, one can show (\ref{eqn: mmr}) reduces to (\ref{eqn: mr}).
Therefore, proof of Remark~\ref{rem: liufinal}(b) is complete.
%Moreover, we prove Remark~\ref{rem: liufinal} under (A1), (A2) and
%(A3). Remark~\ref{remark: truncation} justifies that (A2) can be
%replaced by (A4) and (A5).
%where $K(z,\theta)$ satisfies (\ref{eqn: k}). The pair of equations (
%\ref{eqn: k}) and (\ref{eqn: m}) is the same pair of equations $(2.4)$
%and $(2.5)$ in \citet{LAP2013}. In fact, more generally???? ????
% Now, note that $m(z) = - \sum_{r=1}^{\infty} z^{-r}
%\varphi(d_{uq}^{r-1})$, which is (\ref{eqn: yk}) after replacing $h(
%\lambda, \theta)$ by $1$. All the above calculations remain valid if
%we replace $h(\lambda, \theta)$ by $1$. Hence we have
%
%Therefore, the Stieltjes transformation $m(z)$ of $\frac{1}{2}(\hat{
%\Gamma}_{u}+\hat{\Gamma}_{u}^{*})$ can be expressed as the solution of
%two functional equations (\ref{eqn: k}) and (\ref{eqn: m}).
%Hence the proof of Remark~\ref{rem: liufinal} is complete.
%\vskip5pt

%s5.6 #&#
\subsection{Justification for Example \texorpdfstring{\protect\ref{rem: marchenko}}{1}}
\label
{subsection: marchenko}

%Note that by Corollary~\ref{rem: sympoly}(b), it is enough to work
%under Assumption (A3) or (A4).\vskip5pt
%
(a) Observe that, by (\ref{eqn: lsdmoments}),
\[
\lim p^{-1}E \operatorname{Tr}(\hat{\Gamma}_{0})^{h}
= y^{-1}(1+y)\varphi (\gamma_{00}^{h} ) \qquad\forall
h \geq1,
\]
%
%Observe that, by Theorem~\ref{thmm: lsd}, the limit of $n^{2}(n+p)^{-2}
%\bar{\hat{\Gamma}}_{u}\bar{\hat{\Gamma}}_{u}^{*}$ is $\gamma_{u0}
%\gamma_{u0}^{*}$.
%$$ n\bar{\hat{\Gamma}}_{i}(\varepsilon) = \bar{I}_{p}W\underbar{P}_{i}W
%\bar{I}_{p}.$$
where by (\ref{eqn: gammafree}), % By Theorem~\ref{thmm: lsd}, the
%limit
%of $n(n+p)^{-1}\bar{\hat{\Gamma}}_{0}$ is $\gamma_{00}$, % For the
%sample variance-covariance matrix $\hat{\Gamma}_{0}(\varepsilon)$, we
%have
%$$n\bar{\hat{\Gamma}}_{0}(\varepsilon) = \bar{I_{p}}W\underbar{I}_{n}W
%\bar{I_{p}},$$
%where $I_{p}$ and $I_{n}$ are identity matrices of order $p$ and $n$
%respectively. Therefore,
%by (\ref{eqn: gammafree}),
$\gamma_{00} = \gamma_{00}^{*} = (1+y)a_{0} sc_{0} s a_{0}$. By (\ref
{eqn: limitc}) and (\ref{eqn: limita}), $a_0$ and $c_0$ are both
Bernoulli random variables with success probabilities $y(1+y)^{-1}$ and
$(1+y)^{-1}$, respectively.
%$$\varphi(a_{0}^k) = \frac{y}{1+y} \ \operatorname{and}\ \
%\varphi(c_{0}^k)
%= \frac{1}{1+y} \ \ \operatorname{for all}\ \ k \geq1.$$ %Note that
%(see page
%$144$ of \citet{NS2006}),
%$$\operatorname{$\#$}\ \{\pi\in \mathit{NC}(h):\ \pi\ \operatorname{has $k$
%blocks}\} =
%\frac{1}{k}{h-1 \choose k-1}{h \choose k-1}.$$
Hence, by (\ref{eqn: freemoments}), the $h$th moment of the
LSD of $\hat{\Gamma}_{0}$ is given by
%
%e5.27 #&#
\begin{equation}
\label{eqn: example1a1}  \frac{(1+y)^{h+1}}{y} \sum_{\pi\in \mathit{NC}(h)}
\varphi_{\pi
} \bigl[a_{0}^{2},a_{0}^{2},
\ldots, a_{0}^{2} \bigr]\varphi_{K(\pi
)}[c_{0},c_{0},
\ldots, c_{0}].
\end{equation}
Note that if $\pi\in \mathit{NC}(h)$ has $k$ blocks, then
\begin{eqnarray*}
\varphi_{\pi} \bigl[a_{0}^{2},a_{0}^{2},
\ldots, a_{0}^{2} \bigr] &=& \varphi _{\pi
}[a_{0},a_{0},
\ldots, a_{0}] = y^{k}(1+y)^{-k},
\\
\varphi_{\pi} [c_0, c_0,\ldots,
c_0] &=& (1+y)^{-k}.
\end{eqnarray*}
By (9.18) on page $148$ in \citet{NS2006}, if $\pi\in \mathit{NC}(h)$ has $k$
blocks then $K(\pi)$ has $(h-k+1)$ many blocks. Therefore, %$
%\varphi_{K(
%\pi)}[c_{0},c_{0},\ldots, c_{0}] = (1+y)^{-h-1+k}$. Hence
(\ref{eqn: example1a1}) equals
\begin{eqnarray*}
\sum_{k=1}^{h} \mbox{$\#$} \bigl\{\pi
\in \mathit{NC}(h): \pi \mbox{ has $k$ blocks} \bigr\} y^{k-1} = \sum
_{k=1}^{h}\frac{1}{k}\pmatrix{h-1 \cr k-1}\pmatrix {h
\cr k-1} y^{k-1},
\end{eqnarray*}
which is the $h$th moment of the Mar\v{c}enko--Pastur law [see (\ref
{eqn: MPmoments})]. %The first equality holds by $(9.18)$ on page $148$
%in \citet{NS2006} and f
For the last equality, see page $144$ of \citet{NS2006}. This proves
(a).

(b) Observe that, by (\ref{eqn: lsdmoments}),
\begin{eqnarray*}
\lim p^{-1}E \operatorname{Tr} \bigl( \bigl(np^{-1}
\bigr)^2\hat{\Gamma }_{u}(\varepsilon )\hat{
\Gamma}_{u}^{*}(\varepsilon) \bigr)^{h} =
y^{-(2h+1)}(1+y)\varphi \bigl(\gamma_{u0}\gamma_{u0}^{*}
\bigr)^{h} \qquad\forall h \geq1,
\end{eqnarray*}
%
%Observe that, by Theorem~\ref{thmm: lsd}, the limit of $n^{2}(n+p)^{-2}
%\bar{\hat{\Gamma}}_{u}\bar{\hat{\Gamma}}_{u}^{*}$ is $\gamma_{u0}
%\gamma_{u0}^{*}$.
%$$ n\bar{\hat{\Gamma}}_{i}(\varepsilon) = \bar{I}_{p}W\underbar{P}_{i}W
%\bar{I}_{p}.$$
where, by (\ref{eqn: gammafree}), $\gamma_{u0} = (1+y)
a_{0}sc_{u}sa_{0}$ and $\gamma_{u0}^{*} = (1+y) a_{0}sc_{u}^{*}sa_{0}$.
Since the marginal distribution of all the $c_{i}$'s are same for $i
\geq1$, using free independence of $\operatorname{Span}\{s\}$, $\{
a_{i},a_{i}^{*}\}$ and $\{c_{i},c_{i}^{*}\}$, the LSD of $ (\frac
{n}{p} )^2\hat{\Gamma}_{u}(\varepsilon)\hat{\Gamma
}_{u}^{*}(\varepsilon)$ are identical for all $u \geq1$.

Now we show that the LSD is the free Bessel law. Let
\[
\mathit{NCE}(2n) = \bigl\{\pi\in \mathit{NC}(2n): \mbox{ every block of $\pi$ has even cardinality}
\bigr\}.
\]
%
%Also, by (\ref{eqn: limita}), for all $k \geq1$, $\varphi(a_{0}^{k})
%= \frac{y}{1+y}$. For any monomial $m(\{c_{i},c_{i}^{*}\})$, by (
%\ref{eqn: limitc}), we know $\varphi( m ( \{c_i, c_i^*\}))$.
% from $\operatorname{Span}\{c_{i},c_{i}^{*}\}$, by (\ref{eqn: limitc}),
%\begin{equation} \label{eqn: limitc1} \varphi( m ( \{c_i, c_i^*\})) =
%\left\{
%\begin{array}{ll}
%\frac{1}{1+y} ,& {\rm if there is same number of} c_i {\ \rm and } \
%c_i^* \\
%0, & {\rm otherwise}. \end{array} \right.
%\end{equation}
%Hence using (\ref{eqn: sympoly}) and
By (\ref{eqn: freemoments}), the $h$th order moment of the LSD of
$
(\frac{n}{p} )^{2}\hat{\Gamma}_{u}(\varepsilon)\hat{\Gamma
}_{u}^{*}(\varepsilon)$ is given by
%
%e5.28 #&#
\begin{equation}\qquad
\label{eqn:
example1b1}  \frac{(1+y)^{2h+1}}{y^{2h+1}} \sum_{\pi\in \mathit{NC}(2h)}
\varphi _{K(\pi
)}[a_{0},a_{0},\ldots,
a_{0}]\varphi_{\pi
} \bigl[c_{i},c_{i}^{*},c_{i},c_{i}^{*},
\ldots, c_{i},c_{i}^{*} \bigr].
\end{equation}
%
%& =& y^{-2h}\frac{(1+y)^{2h+1}}{y} \sum_{\pi\in \mathit{NCE}(2h)} \varphi_{
%\pi}[a_{0},a_{0},\ldots, a_{0}]\varphi_{K(
%\pi)}[c_{i},c_{i}^{*},c_{i},c_{i}^{*}\ldots, c_{i},c_{i}^{*}]\\
Note that $\varphi_{\pi}[c_{i},c_{i}^{*},c_{i},c_{i}^{*},\ldots,
c_{i},c_{i}^{*}] = 0$ if $\pi\in \mathit{NC}(2h)-\mathit{NCE}(2h)$. If $\pi\in
\mathit{NCE}(2h)$ has $k$ many blocks, then $\varphi_{\pi
}[c_{i},c_{i}^{*},c_{i},c_{i}^{*},\ldots, c_{i},c_{i}^{*}] = (1+y)^{k}$.
Note that by $(9.18)$ on page $148$ in \citet{NS2006}, $K(\pi)$ has
$2h+1-k$ blocks, and hence $\varphi_{K(\pi)}[a_{0},a_{0},\ldots, a_{0}]
= y^{2h+1-k}(1+y)^{2h+1-k}$. Therefore, (\ref{eqn: example1b1}) equals
%
%e5.29 #&#
\begin{eqnarray}
 &&y^{-2h}\sum_{k=1}^{h}
\mbox{$\#$} \bigl\{\pi\in \mathit{NCE}(2h): \pi \mbox { has $k$ blocks} \bigr\}
y^{2h+1-k-1}\nonumber
\\
\eqntext{ \mbox{[by (\ref{eqn: limitc}) and (\ref{eqn: limita})]}}
\\
&&\qquad =  \sum_{k=1}^{h}\frac{1}{k}{h-1
\choose k-1} {2h \choose k-1} y^{-k},\nonumber
\end{eqnarray}
where
%as in Corollary~\ref{rem: marchenko}(a), the first equality holds by
%$(9.18)$ on page $148$ in \citet{NS2006} and
the last equality follows from Lemma $4.1$ of \citet{E1980}.
%$$\mbox{$\#$} \{\pi\in \mathit{NCE}(2h):\ \pi\ \mbox{has $k$ blocks}\} =
%\frac{1}{k}{h-1 \choose k-1}{2h \choose k-1}.$$
The final expression is indeed the $h$th moment of the free
Bessel$(2,y^{-1})$ law.
%If $p \to\infty$ such that $y=0$, then %following the arguments
%leading to Remark~\ref{rem: bassel},
%the $h$-th order moment of the LSD of $\left(\frac{n}{p}\right)\hat{
%\Gamma}_{i}(\varepsilon)\hat{\Gamma}_{i}^{*}(\varepsilon)$ is $
%\frac{1}{h}{2h \choose h-1}$. Hence, the LSD is the Mar
%\v{c}enko-Pastur law with parameter $1$.
This proves (b) under (A2).
%The third equality follows from (\ref{eqn: limitc1}).
%\end{remark}
%\vskip10pt
%\begin{remark} \label{rem: y0bessel}

%\end{remark}
%\begin{remark}

%\begin{remark}
%\vskip10pt
%

(c) By (\ref{eqn: k}) and (\ref{eqn: m}), the
Stieltjes transform of LSD of $\frac{1}{2}(\hat{\Gamma}_u + \hat
{\Gamma
}^{*}_u)$ satisfies,
%
%e5.30 #&#
\begin{equation}
\label{eqn: stma0gamma1} m(z) = - \biggl(z - \frac{1}{2\pi}\int_{0}^{2\pi}
\frac{\cos\theta
\,d\theta}{1+ym(z)\cos\theta} \biggr)^{-1}.
\end{equation}
Now by contour integration, it can be shown that
\[
\frac{1}{2\pi}\int_{0}^{2\pi}\frac{\cos\theta \,d\theta
}{1+ym(z)\cos
\theta}
%& = & \frac{1}{2\pi ym(z)}\int_{0}^{2\pi}\frac{(1+ym(z)\cos
%\theta)-1}{1+ym(z)\cos\theta} d\theta\nonumber\\
%& = & \frac{1}{ym(z)} - \frac{1}{2\pi ym(z)}\int_{0}^{2\pi}
%\frac{1}{1+ym(z)\cos\theta} d\theta\nonumber\\
%& = & \frac{1}{ym(z)} - \frac{1}{2\pi ym(z)}\int_{0}^{2\pi}
%\frac{1}{1+(1/2)ym(z)(e^{i\theta}+e^{-i\theta})} d\theta\nonumber\\
%& = & \frac{1}{ym(z)} - \frac{1}{2\pi i ym(z)}\int_{\gamma}
%\frac{1}{1+(1/2)ym(z)(\omega+ \omega^{-1})}\frac{d\omega}{\omega}
%\nonumber\\
%& = & \frac{1}{ym(z)} - \frac{2}{2\pi i y^{2}m^{2}(z)}\int_{\gamma}
%\frac{d\omega}{\omega^{2}+2(ym(z))^{-1}\omega+1} \nonumber\\
=  \frac{1}{ym(z)} - \frac{2}{y^{2}m^{2}(z)} \frac{1}{\omega_1 -
\omega_2}, %\frac{yG(z)}{2\sqrt{1-y^{2}G^{2}(z)}}
\nonumber
%& = & -\frac{1}{yG(z)} - \frac{1}{yG(z)}\frac{1}{
%\sqrt{1-y^{2}G^{2}(z)}}. \nonumber
\]
where $\omega_1$ and $\omega_2$ are two roots of $\omega
^{2}+2(ym(z))^{-1}\omega+1=0$ with $|\omega_1|>1$, $|\omega_2|<1$ and
$(\omega_1-\omega_2)^{-2} = \frac{y^{2}m^{2}(z)}{4(1-y^{2}m^{2}(z))}$.
%%Also $\gamma= \{e^{i\theta}: 0\leq\theta< 2\pi\}$ i.e. contour of
%unit circle.
Therefore, by (\ref{eqn: stma0gamma1}), we have
\begin{eqnarray*}
\label{eqn: prec1}  \frac{-1}{m(z)} &=& z -\frac{1}{ym(z)} +
\frac
{2(\omega_1 - \omega_2)^{-1}}{y^{2}m^{2}(z)} %\frac{1}{} %
%\frac{1}{yG(z)}\frac{1}{\sqrt{1-y^{2}G^{2}(z)}}
%\nonumber\\
%&\implies& -\frac{1}{m(z)}\left(1-\frac{1}{y}
%\right)-z = \frac{2}{y^{2}m^{2}(z)} \frac{1}{\omega_1 - \omega_2}%
%\frac{1}{yG(z)}\frac{1}{\sqrt{1-y^{2}G^{2}(z)}}
% \nonumber\\
% &\implies& -y\left(1-\frac{1}{y}\right)-zym(z) =
%\frac{2}{ym(z)} \frac{1}{\omega_1 - \omega_2}%\frac{1}{
%\sqrt{1-y^{2}G^{2}(z)}}
%\nonumber\\
%&\implies& (1-y)-zym(z) = \frac{2}{ym(z)} \frac{1}{
%\omega_1 - \omega_2}%\frac{1}{\sqrt{1-y^{2}G^{2}(z)}}
%\nonumber\\
\\
&\implies&  \bigl((1-y)-zym(z) \bigr)^{2} \bigl(1-y^{2}m^{2}(z)
\bigr) = 1.
%\mbox{or,} 1 & = & \left((1-y)^2 + z^{2}y^{2}G^{2}(z)
%+ 2(1-y)yzG(z)\right)\left(1-y^{2}G^{2}(z) \right) \nonumber\\
%\mbox{or,} - z^{2}y^{4}G^{4}(z) -
%2y^{3}(1-y)zG^{3}(z)- ((1-y)^2 - z^2) y^{2}G^{2}(z) + 2y(1-y)zG(z) +
%(1-y)^{2}-1 = 0 \nonumber
\end{eqnarray*}

Hence, Example~\ref{rem: marchenko}(c) follows. % Theorem~2.1 of \citet{LAP2013} also agrees with (\ref{eqn: bai}). %
%\citet{JWBN2014}.} % (see (\ref{eqn: bai}).}

%s5.7 #&#
\subsection{Justification for Example \texorpdfstring{\protect\ref{rem: consta}}{3}}
\label
{subsection: constaproof} By Remark~\ref{remark: truncation}(a) it is
enough to work under assumption (A2). Note that the LSD of $\psi_j$
is $\delta_{\lambda_j}$. % (degenerated at $\lambda_j$). For any $i$,
%let $P_i$
%be the $n \times n$ matrix with upper $i$th diagonal $1$. Note
%that $P_0=I_n$.
We can write %????factor of $n$ does nto seem to be ok??
\[
n{\Delta}_{u} = Z \Biggl(\sum_{j=0}^{q}
\sum_{j^\prime=0}^{q}\lambda _{j}
\lambda_{j^\prime}P_{j-j^\prime+u} \Biggr)Z^*, \qquad n{\Delta}_{u}^{*}
= Z \Biggl(\sum_{j=0}^{q}\sum
_{j^\prime
=0}^{q}\lambda _{j}
\lambda_{j^\prime}P_{j-j^\prime+u}^{*} \Biggr)Z^*.
\]
By Lemma $7.1$ of the supplementary file \citet{BB2014freesup}, LSD of
$\frac{1}{2}(\hat{\Gamma}_{u}+\hat{\Gamma}_{u}^{*})$ and $\frac
{1}{2}({\Delta}_{u}+{\Delta}_{u}^{*})$ are identical.
Moreover, %By Remark~\ref{rem: sympoly}, the LSD of $
%\frac{1}{2}(\bar{\hat{\Gamma}}_{i}+\bar{\hat{\Gamma}}_{i}^{*})$ is
%given by
%
\[
\frac{1}{2} \bigl({{\Delta}}_{u}+{{\Delta}}_{u}^{*}
\bigr) = n^{-1}Z \Biggl(\frac
{1}{2}\sum
_{j=0}^{q}\sum_{j^\prime=0}^{q}
\lambda_{j}\lambda _{j^\prime
} \bigl(P_{j-j^\prime+u}+P_{j-j^\prime+u}^{*}
\bigr) \Biggr)Z^*,
\]
whose LSD is a compound free Poisson [see discussions around
(\ref{eqn: zaz})] and by~(\ref{eqn: zaz}), its $r$th order free
cumulant %of LSD of $\frac{1}{2}({\hat{\Gamma}}_{i}+{\hat{
%\Gamma}}_{i}^{*})$
is given by
\begin{eqnarray*}
%\label{eqn: constacum}
%K_{ir} =
&&
y^{r-1}\lim\frac{1}{n}\Tr
\Biggl(\frac{1}{2}\sum_{j,j^\prime
=0}^{q}
\lambda_{j}\lambda_{j^\prime} \bigl(P_{j-j^\prime+u}+P_{j-j^\prime
+u}^{*}
\bigr) \Biggr)^{r}\\
&&\qquad = y^{r-1}E_{\theta} \bigl(\cos(u
\theta)\tilde {h}(\lambda, \theta) \bigr)^r,
\end{eqnarray*}
%
%with $f_{j}(\lambda) = \lambda_{j}\ \forall j$ and
where $\tilde{h}$ is as given in (\ref{eqn: h}) and $\theta\sim
U(0,2\pi)$. Hence, Example~\ref{rem: consta} is justified.

%s5.8 #&#
\subsection{Proof of Remark \texorpdfstring{\protect\ref{rem: paffel}}{3.3}} \label{subsection: paffel}
It is enough to show that %Now, we will show that
(\ref{eqn: constacum}) justifies Theorem~$2.1$ in \citet{LAP2013} when
$\psi_j = \lambda_j$.

For a variable $a$, define the free cumulant generating function
%
%e5.31 #&#
\begin{equation}
\label{eqn: cgf} C(z) = 1+\sum_{r=1}^{\infty}k_{r}
z^{r},
\end{equation}
where $k_{r}$ is the $r$th order free cumulant of $a$.
%where $k_{ur}$ is as in (\ref{eqn: constacum}).
One can show that
%
%e5.32 #&#
\begin{equation}
\label{eqn: cm} -C \bigl(-m(z) \bigr) = zm(z),
\end{equation}
[see $(12.5)$ on page $198$ of \citet{NS2006}] where $m(z)$ is the
Stieltjes transformation of $a$. % $\frac{1}{2}(\hat{\Gamma}_{u} + \hat{
%\Gamma}_{u}^{*})$.

Here, we consider $a$ to be distributed as LSD of $\frac{1}{2}(\hat
{\Gamma}_{u} + \hat{\Gamma}_{u}^{*})$. Our goal is to find the power
series expansion of $C(-m(z))$ from $m(z)$ in (\ref{eqn: mr}) and show
that the coefficient of $(-m(z))^{r}$ in that expansion is same as
$k_{ur}$ in (\ref{eqn: constacum}).

As $\psi_j = \lambda_jI_p$, (\ref{eqn: hr}) reduces to (\ref{eqn: h}),
and hence (\ref{eqn: mr}) and (\ref{eqn: kr}) reduce to
%
%e5.33 #&#
%e5.34 #&#
\begin{eqnarray}
\label{eqn: mmmr} \frac{1}{m(z)} + z %\nonumber\\
&=& E_{\theta} \biggl(
\frac{\cos(u\theta)\tilde{h}(\lambda,\theta
) }{1+
y \cos(u \theta) K(z,\theta)} \biggr)\qquad\mbox{where}
\\
%\nonumber
\label{eqn: kkkr} K(z,\theta) %\nonumber\\
%&=& \bigg( E_{\theta^\prime}\big(\frac{\cos(u\theta^\prime) h_{1}(
%\alpha,\theta^\prime) }{1+ y \cos(u \theta^\prime) K(z,\theta^\prime)}
%\big)-z\bigg)^{-1} h_{1}(\alpha, \theta)
&=&
m(z)\tilde{h}(\lambda, \theta). %\nonumber\\
% h_{1}(\alpha, \theta) = |\sum_{j=0}^{q} e^{ij\theta}
%\lambda_{j}|^{2}.
\end{eqnarray}
Hence, by (\ref{eqn: cm}),
\begin{eqnarray}\label{eqn: cumexp}
C \bigl(-m(z) \bigr) &=& - zm(z) = -m(z) \biggl( \frac{1}{m(z)} + z \biggr)
+ 1
\nonumber
\\
&=& - E_{\theta} \biggl(\frac{\cos(u\theta) m(z) \tilde{h}(\lambda,
\theta)}{1+ y \cos(u \theta) K(z,\theta)} \biggr) + 1 \qquad\mbox{by
(\ref
{eqn: mmmr})}
\nonumber
\\[-8pt]
\\[-8pt]
\nonumber
&=& - E_{\theta} \biggl(\frac{\cos(u\theta) m(z) \tilde{h}(\lambda,
\theta) }{1+ y \cos(u \theta) m(z) \tilde{h}(\lambda, \theta
)} \biggr) + 1\qquad \mbox{by
(\ref{eqn: kkkr})}
\nonumber
\\
%&=& -\frac{1}{2\pi}\int_{0}^{2\pi}\frac{\cos(u\theta)m(z)h(\lambda,
%\theta)d\theta}{1+y\cos(u\theta)m(z)h(\lambda,\theta)} + 1 \\
 &=& 1+\sum_{r=1}^{\infty}y^{r-1}E_{\theta}
\bigl(\cos (u\theta)\tilde{h}(\lambda, \theta) \bigr)^{r} \bigl(-m(z)
\bigr)^{r}.\nonumber
\end{eqnarray}

%\vskip5pt By (\ref{eqn: k}) and (\ref{eqn: m}), we have
%\begin{equation} \label{eqn: simplifystieltjes}
%\frac{1}{2\pi}\int_{0}^{2\pi}\frac{\cos(u\theta)K(z,\theta)d\theta}{1+y
%\cos(u\theta)K(z,\theta)} -1 = zm(z) = -C(-m(z)).
%\end{equation}
% Moreover, since $\psi_{j} = \lambda_{j}I_{p}$ for all $j$,
%we have $K(z,\theta) = m(z)h(\lambda,\theta)$, where $h(\lambda,
%\theta)$ is given by (\ref{eqn: h}). Hence, by (\ref{eqn:
%simplifystieltjes}), we have
%\begin{eqnarray}
%\nonumber C(-m(z)) &=& -\frac{1}{2\pi}\int_{0}^{2\pi}\frac{\cos(u
%\theta)m(z)h(\lambda,\theta)d\theta}{1+y\cos(u\theta)m(z)h(\lambda,
%\theta)} + 1 \\
%\label{eqn: cumexp} &=& 1+\sum_{r=1}^{\infty}y^{r-1}E_{\theta}(\cos(u
%\theta)h(\lambda,\theta))^{r}(-m(z))^{r},
%\end{eqnarray}
%where $C(\cdot)$ is the free cumulant generating function
%defined in (\ref{eqn: cgf}).
Therefore, for all $r\geq1$, the coefficient of $(-m(z))^{r}$ in (\ref
{eqn: cumexp}) agrees with $k_{ur}$ in (\ref{eqn: constacum}). This
completes the proof.

%s5.9 #&#
\subsection{Detailed calculation for Example \texorpdfstring{\protect\ref{example:
5}}{5} in
Section \texorpdfstring{\protect\ref{subsubsec: trace}}{4.2.3}} \label{subsec: detailtrace}
Note that $X_t = \varepsilon_t + \varepsilon_{t-1}$, where
$\varepsilon_t \sim\mathcal{N}(0,I_p)$ and $I_p$ is the identity
matrix of order $p$. Suppose $n=p$. Therefore, for every $i \geq1$,
$(\varepsilon_{t,i} +\varepsilon_{t-1,+i}) \sim\mathcal{N}(0,2)$,
and hence
%
%e5.35 #&#
%e5.36 #&#
\begin{eqnarray}\qquad
\label{eqn: mm4} E(\varepsilon_{t,i} +\varepsilon_{t-1,+i})^{2}
&=& 2, \qquad E(\varepsilon_{t,i} +\varepsilon_{t-1,+i})^{4}
= 12,
\\
\label{eqn: ma1ex} E \bigl(\operatorname{Tr}(\hat{\Gamma}_{0}) \bigr) &=&
n^{-1} \sum_{t_1,i_1}E \bigl(X_{t_1,i_1}^{2}
\bigr) = n^{-1} \sum_{t_1,i_1}E(\varepsilon
_{t_1,i_1} + \varepsilon_{t_1-1, i_1})^{2}
\nonumber
\\[-8pt]
\\[-8pt]
\nonumber
&=& 2(n-1)
\end{eqnarray}
and
\begin{eqnarray*}
&& E \bigl(\operatorname{Tr}(\hat{\Gamma}_{0}) - E
\operatorname {Tr}(\hat{\Gamma}_{0}) \bigr)^{2}
\\
&&\qquad= E \bigl(\operatorname{Tr}(\hat{\Gamma}_{0}) \bigr)^{2}
- \bigl(E\operatorname {Tr}(\hat{\Gamma}_{0}) \bigr)^{2}
\\
&&\qquad= n^{-2} \sum_{t_1,t_2,i_1,i_2} E
\bigl(X_{t_1,i_1}^{2} X_{t_2,i_2}^{2} \bigr) -
4(n-1)^{2} \qquad\bigl[\mbox{by (\ref{eqn: ma1ex})}\bigr]
\\
&&\qquad= n^{-2} \sum_{t_1 \neq t_2,i_1, i_2} E
\bigl(X_{t_1,i_1}^{2} \bigr)E \bigl( X_{t_2,i_2}^{2}
\bigr) + n^{-2} \sum_{t_1=t_2,i_1=i_2} E
\bigl(X_{t_1,i_1}^{4} \bigr)
\\
&&\quad\qquad{} + n^{-2} \sum_{t_1 = t_2,i_1\neq i_2} E
\bigl(X_{t_1,i_1}^{2} \bigr)E \bigl( X_{t_2,i_2}^{2}
\bigr) - 4(n-1)^{2}
\\
%&= & 4n(n-1) + 12 + 4(n-1) - 4n^{2} , \mbox{(by (\ref{eqn: mm2}) and (
%\ref{eqn: mm4}))} \nonumber\\
&&\qquad= 4(n-1) (n-2) + 12 + 4\frac{(n-1)^{2}}{n}
-4(n-1)^{2}\qquad\bigl[\mbox{by (\ref{eqn: mm4})}\bigr]
\\
&&\qquad= 12 + 4\frac{n-1}{n} \bigl(n^2 -2n +n-1-n^2 +n
\bigr) = 12 - 4\frac{n-1}{n}
\\
&&\qquad\to 8.
\end{eqnarray*}

\section*{Acknowledgements} We thank Octavio Arizmendi and
Carlos Vargas Obieta for interesting and helpful discussions on free
probability. We thank Debashis Paul for interesting discussions on
their article. We thank all the three referees for their extremely
constructive comments which have led to a very significant improvement
in the article in terms of both substance and presentation.
%bringing up some important points.

\begin{supplement}[id=suppA]
\stitle{Supplement to ``Large sample behaviour of high dimensional autocovariance matrices''}
\slink[doi]{10.1214/15-AOS1378SUPP} %[doi,text={...}] - jei reikia
%suskaldyti doi
\sdatatype{.pdf}
\sfilename{aos1378\_supp.pdf}
\sdescription{In this supplement, we provide additional technical
details and simulations.}
\end{supplement}

%\begin{supplement}[id-suppA]
%\stitle{Supplement to ``Large sample behaviour of high dimensional
%autocovariance matrices"}
%\slink[doi]{COMPLETED BY THE TYPESETTER}
%\sdescription{In this supplement, we provide additional technical
%details and simulations.}
%\end{supplement}
%
%\bibliographystyle{imsart-nameyear}
%\bibliography{reference}
%
% imsref loaded by akundreckaite, 2015-10-21 09:45:08
%

%\begin{appendix}
%\section{}
%\end{appendix}

%\section*{Acknowledgments}

%\begin{thebibliography}{99}
%\bibitem{r1}
%\bibitem{r1}
%\end{thebibliography}

\printaddresses
\end{document}